\documentclass[12pt]{article}

\usepackage[left=1in,top=1in,right=1in,bottom=1in]{geometry}

\usepackage{graphicx}
\usepackage{amsmath}
\usepackage{amssymb}
\usepackage{amsthm}
\usepackage{latexsym}
\usepackage{bm}
\usepackage{color}
\usepackage{fancyhdr}
\usepackage{stmaryrd}
\usepackage{mathrsfs}
\usepackage{url}
\usepackage{array}
\usepackage{algorithm}
\usepackage{algpseudocode}
\usepackage{multirow}

\SetSymbolFont{stmry}{bold}{U}{stmry}{m}{n}
\newtheorem{remark}{Remark}

\numberwithin{equation}{section}

\newcolumntype{P}[1]{>{\centering\arraybackslash}p{#1}}

\begin{document}
\title{A robust and efficient iterative method for hyper-elastodynamics with nested block preconditioning}
\author{Ju Liu and Alison L. Marsden\\
\textit{\small Department of Pediatrics (Cardiology), Bioengineering, and}\\
\textit{\small Institute for Computational \& Mathematical Engineering, Stanford University}\\
\textit{\small Clark Center E1.3, 318 Campus Drive, Stanford, CA 94305, USA}\\
\small \textit{E-mail address:} liuju@stanford.edu, amarsden@stanford.edu
}

\date{}
%
%
\maketitle

\section*{Abstract}
We develop a robust and efficient iterative method for hyper-elastodynamics based on a novel continuum formulation recently developed in \cite{Liu2018}. The numerical scheme is constructed based on the variational multiscale formulation and the generalized-$\alpha$ method. Within the nonlinear solution procedure, a block factorization is performed for the consistent tangent matrix to decouple the kinematics from the balance laws. Within the linear solution procedure, another block factorization is performed to decouple the mass balance equation from the linear momentum balance equations. A nested block preconditioning technique is proposed to combine the Schur complement reduction approach with the fully coupled approach. This preconditioning technique, together with the Krylov subspace method, constitutes a novel iterative method for solving hyper-elastodynamics. We demonstrate the efficacy of the proposed preconditioning technique by comparing with the SIMPLE preconditioner and the one-level domain decomposition preconditioner. Two representative examples are studied: the compression of an isotropic hyperelastic cube and the tensile test of a fully-incompressible anisotropic hyperelastic arterial wall model. The robustness with respect to material properties and the parallel performance of the preconditioner are examined. \\

\noindent \textbf{Keywords:} Variational multiscale method, Saddle-point problem, Nested iterative method, Block preconditioner, Anisotropic hyperelasticity, Arterial wall model 




\section{Introduction}
\label{sec:introduction}
In our recent work \cite{Liu2018}, a unified continuum modeling framework was developed. In this framework, hyperelastic solids and viscous fluids are distinguished only through the deviatoric part of the Cauchy stress, in contrast to prior modeling approaches. In our derivation, the Gibbs free energy, rather than the Helmholtz free energy, is chosen as the thermodynamic potential, resulting in a unified model for compressible and incompressible materials. A beneficial outcome of the modeling framework is that it naturally allows one to apply a computational fluid dynamics (CFD) algorithm to solid dynamics, or vice versa. In our work \cite{Liu2018}, the variational multiscale (VMS) analysis, a mature numerical modeling approach in CFD \cite{Hughes1995}, is taken to design the spatial discretization for solid dynamics. This numerical model provides a stabilization mechanism that circumvents the Ladyzhenskaya-Babu\v{s}ka-Brezzi (LBB) condition for equal-order interpolations. In particular, it allows one to use low-order tetrahedral elements, even for fully incompressible materials. This gives us the maximum flexibility in geometrical modeling and mesh generation. 

In this work, we build upon the proposed unified formulation to develop a robust and efficient iterative method. Traditional black-box preconditioners are non-robust, and the convergence rate of the linear solver drops significantly under certain conditions. The lack of robustness may be attributed to the saddle-point nature of the problem. Algebraic preconditioners built based on incomplete factorizations are prone to fail due to zero-pivoting; one-level domain decomposition preconditioners do not perform well due to its locality. In this work, we design a preconditioning technique tailored for the VMS formulation for hyper-elastodynamics \cite{Liu2018}. The design of the preconditioner is based on a nested block factorization of the consistent tangent matrix in the Newton-Raphson iteration.  A block factorization is performed in the nonlinear solution procedure to decouple the kinematics from the balance laws \cite{Scovazzi2016}. The resulting $2\times 2$ block matrix is further factorized in the linear solution procedure.  This strategy is, in part, related to the classical projection method \cite{Chorin1968,Tem1969} and the block preconditioning technique \cite{Benzi2005,Elman2014,Turek1999} that have been widely used in the CFD community. We examine the solver performance for both isotropic and anisotropic hyperelastic models. The significance of this work is that it paves the way towards robust, efficient, and scalable implicit solver technology for biomechanics and monolithic fluid-solid interaction (FSI) simulations \cite{Liu2018}. In the rest part of this section, we give an overview of the background and an outline of the work.

\subsection{Projection method and block preconditioners}
The development of efficient solver techniques for multiphysics problems has been an active area of research in recent years \cite{D.E.Keyes2013}. One simple but important prototype multiphysics problem is the Stokes or the Navier-Stokes equations, representing the coupling between the mass conservation and the balance of the linear momentum for incompressible flows. In the late 1960s, the Chorin-Teman projection method \cite{Chorin1968,Tem1969} was proposed to solve for the pressure and the velocity separately based on the Helmholtz decomposition. Since then, the projection method and its variants have attracted concentrated research and lead to a voluminous literature \cite{Kim1985,Kan1986,Karniadakis1991,Guermond2003}. The projection method is attractive because the nonlinear system of equations is decomposed into a series of linear elliptic equations. Although this method has attracted significant attention, it still poses several major challenges. One critical issue is that the physics-based splitting necessitates the introduction of an artificial boundary condition for the pressure. There is no general theory to guide the choice of the artificial boundary conditions, and most likely this artificial boundary condition limits the solution accuracy. For an overview of the projection method, the readers are referred to the review article \cite{Guermond2006}.

In recent years, it has been realized that one can invoke an arbitrary time stepping scheme (e.g. fully implicit) and achieve the decoupling of physics within the linear solver. Indeed, in each iteration of the Krylov subspace method, one only needs to solve with a preconditioner and perform a matrix-vector multiplication to construct the new search direction. Therefore, if the preconditioner is endowed with a block structure, one may sequentially solve each block matrix with less cost. It has been pointed out that the Chorin-Teman projection method is closely related to a block preconditioner \cite{Perot1993}. Consider a matrix problem with a $2\times 2$ block structure,
\begin{align*}
\mathcal A :=
\begin{bmatrix}
\boldsymbol{\mathrm A} & \boldsymbol{\mathrm B} \\[0.3mm]
\boldsymbol{\mathrm C} & \boldsymbol{\mathrm D}
\end{bmatrix}.
\end{align*}
This matrix can be factored into lower triangular, diagonal, and upper triangular matrices as follows,
\begin{align*}
\mathcal A = \mathcal L \mathcal D \mathcal U =
\begin{bmatrix}
\boldsymbol{\mathrm I} & \boldsymbol{\mathrm O} \\[0.3em]
\boldsymbol{\mathrm C} \boldsymbol{\mathrm A}^{-1} & \boldsymbol{\mathrm I}
\end{bmatrix}
\begin{bmatrix}
\boldsymbol{\mathrm A} & \boldsymbol{\mathrm O} \\[0.3em]
\boldsymbol{\mathrm O} & \boldsymbol{\mathrm S}
\end{bmatrix}
\begin{bmatrix}
\boldsymbol{\mathrm I} & \boldsymbol{\mathrm A}^{-1} \boldsymbol{\mathrm B} \\[0.3em]
\boldsymbol{\mathrm O} & \boldsymbol{\mathrm I}
\end{bmatrix}.
\end{align*}
The diagonal matrix $\mathcal D$ contains a Schur complement $\boldsymbol{\mathrm S} := \boldsymbol{\mathrm D} - \boldsymbol{\mathrm C} \boldsymbol{\mathrm A}^{-1} \boldsymbol{\mathrm B}$,which acts as an algebraic analogue of the Laplacian operator for the pressure field \cite{Quarteroni2000}. To construct a preconditioner for $\mathcal A$, one needs to provide approximations for $\boldsymbol{\mathrm A}$ and $\boldsymbol{\mathrm S}$ that can be conveniently solved with. The new formulation for hyper-elastodynamics we consider here is similar to the generalized Stokes equations, in which the operator $\boldsymbol{\mathrm A}$ arises from the discretization of a combination of zeroth order and second order differential operators. Thus, $\boldsymbol{\mathrm A}$ is amenable for approximation by a standard preconditioning technique. Due to the presence of $\boldsymbol{\mathrm A}^{-1}$, $\boldsymbol{\mathrm S}$ is a dense matrix. When the matrix $\boldsymbol{\mathrm A}$ represents a discretization of a zeroth-order differential operator, an effective choice is to replace $\boldsymbol{\mathrm S}$ by $\hat{\boldsymbol{\mathrm S}} := \boldsymbol{\mathrm D} - \boldsymbol{\mathrm C} \left(\textup{diag}\left(\boldsymbol{\mathrm A}\right)\right)^{-1} \boldsymbol{\mathrm B}$ to construct the preconditioner for $\mathcal A$. This choice is closely related to the SIMPLE scheme commonly used in CFD \cite{Elman2008,Patankar1983}. When the matrix $\boldsymbol{\mathrm A}$ represents a discretization of a second-order differential operator, a scaled mass matrix is often effective \cite{Silvester1994}. For more complicated problems, designing a spectrally equivalent preconditioner for the Schur complement is challenging and, in a broad sense, remains an open question. In recent years, progress has been made for problems where $\boldsymbol{\mathrm A}$ is dominated by a discrete convection operator. Notable 
examples include the BFBT preconditioner \cite{Elm1999}, the pressure convection diffusion preconditioner \cite{Kay2002}, and the least squares commutator (LSC) preconditioner \cite{Elman2006}. Based on the Sherman-Morrison formula, a different preconditioner for the Schur complement can be designed for problems with significant contributions from the boundary conditions \cite{Moghadam2013a}. In all, the block preconditioner, as an algebraic interpretation of the projection method, has become increasingly popular, since it does not necessitate ad hoc pressure boundary conditions and allows fully implicit time stepping schemes.

If one can solve the sub-matrices $\boldsymbol{\mathrm A}$ and $\boldsymbol{\mathrm S}$ to a prescribed tolerance, the matrix $\mathcal A$ is solved in one pass without generating a Krylov subspace. This is commonly known as the \textit{Schur complement reduction} (SCR) or \textit{segregated} approach \cite{Benzi2005,May2008,Lun2010,Furuichi2011}. In contrast, the aforementioned strategy, where $\mathcal A$ is solved by a preconditioned iterative method, is referred to as the \textit{coupled} approach \cite{Benzi2005}. For many problems, it is impractical to explicitly construct the Schur complement. Still, the action of the Schur complement on a vector can be obtained in a ``matrix-free" manner (see Algorithm \ref{algorithm:matrix_free_mat_vec_for_S} in Section \ref{subsec:block_pc}). Thus, one can still solve with the Schur complement by iterative methods. To achieve high accuracy, a sufficient number of bases of the Krylov subspace for $\boldsymbol{\mathrm S}$ need to be generated, and this procedure can be prohibitively expensive.

\subsection{Nested preconditioning technique}
The difference between the coupled approach and the segregated approach can be viewed as follows. In the coupled approach, $\boldsymbol{\mathrm S}$ is replaced by a sparse approximation to generate a preconditioner for $\mathcal A$. In the segregated approach or SCR, one strives to solve directly with $\boldsymbol{\mathrm S}$. The distinction between the two approaches is blurred by using the SCR procedure as a preconditioner. In doing so, one does not need to solve with $\boldsymbol{\mathrm S}$ to a high precision, thus alleviating the computational burden. In comparison with the coupled approach, the information of the Schur complement is maintained in the preconditioner (up to the tolerances of SCR), and this will improve the robustness. Therefore, in solving with $\mathcal A$, there are three nested levels. In the outer level, a Krylov subspace method is applied for $\mathcal A$ with a block preconditioner. In the intermediate level, the block preconditioner is applied by solving with the matrices $\boldsymbol{\mathrm A}$ and $\boldsymbol{\mathrm S}$. In the inner level, a solver of $\boldsymbol{\mathrm A}$ is invoked to approximate the action of $\boldsymbol{\mathrm S}$ on a vector. Two mechanisms guarantee and accelerate the convergence. In the outer level, the Krylov subspace method for $\mathcal A$ minimizes the residual of the coupled problem. In the intermediate and inner levels, the SCR procedure is utilized as the preconditioner, which itself can be viewed as an inaccurate solver for $\mathcal A$.

Using SCR as a preconditioner was first proposed within a Richardson iteration scheme \cite{Bank1990}. Due to the symmetry property of that problem, a conjugate gradient method is applied to solve the Schur complement equation. Later, the nested iterative scheme was investigated for CFD problems \cite{Baggag2004,Manguoglu2008}, and the reported results indicate that using the SCR procedure as a preconditioner in a Richardson iteration outperforms the coupled approach with a Krylov subspace method. The nested algorithm was then further investigated using the biconjugate gradient stabilized method (BiCGStab) as the outer solver \cite{Manguoglu2009}. The nested iterative scheme in \cite{Manguoglu2009} uses rather crude stopping criteria for the intermediate and inner solvers. Still, its performance is superior to that of BiCGStab preconditioned by a BFBT preconditioner.

Our investigation of the VMS formulation for hyper-elastodynamics starts with a SIMPLE-type block preconditioner using our in-house code \cite{Moghadam2015}. As will be shown in Section \ref{sec:numerical_results}, the Krylov subspace method with a block preconditioner like SIMPLE is not always robust. This can be attributed to the ignorance of the off-diagonal entries in $\boldsymbol{\mathrm A}$. Because of that, it is appealing to consider preconditioners like LSC, since the off-diagonal information of $\boldsymbol{\mathrm A}$ is maintained. However, non-convergence has been reported for LSC when solving the Navier-Stokes equations with stabilized finite element schemes \cite{Cyr2012}. We then ruled out this option since our VMS formulation involves a similar pressure stabilization term. Consequently, we consider using SCR with relaxed tolerances as a preconditioner. In doing so, the Schur complement is approximated through using an inner solver. In contrast to the nested iterative approaches introduced above, we adopt the following techniques in our study: (1) we use GMRES \cite{Saad1986} and its variant \cite{Saad1993} as the Krylov subspace method in all three levels to leverage their robustness in handling non-symmetric matrix problems; (2) we apply the algebraic multigrid (AMG) preconditioner \cite{Yang2002} for problems at the intermediate level to enhance the robustness of the overall algorithm; (3) we use the sparse approximation $\hat{\boldsymbol{\mathrm S}}$ as a preconditioner when solving with $\boldsymbol{\mathrm S}$. We demonstrate application of this method to hyper-elastodynamics, however we anticipate its general use in CFD and FSI problems in future work.

\subsection{Structure and content of the paper}
The remainder of the study is organized as follows. In Section \ref{sec:elastodynamics}, we state the governing equations of hyper-elastodynamics \cite{Liu2018}. In Section \ref{sec:numerical_formulation}, the numerical scheme is presented. A block factorization for the consistent tangent matrix is performed to reduce the size of the linear algebra problem. In Section \ref{sec:linear_solver}, the nested block preconditioning technique is discussed in detail. In Section \ref{sec:numerical_results}, we present two representative examples to demonstrate the efficacy of the proposed solver technology. The first example is the compression of an isotropic elastic cube \cite{Reese2000}, and the second is the tensile test of a fully incompressible anisotropic hyperelastic arterial wall model \cite{Gasser2006}. Comparisons with other preconditioners are made. We draw conclusions in Section \ref{sec:conclusion}.

\section{Hyper-elastodynamics}
\label{sec:elastodynamics}
In this section, we state the initial-boundary value problem for hyper-elastodynamics, following the derivation in \cite{Liu2018}. Let $\Omega_{\bm X}$ and $\Omega_{\bm x}$ be bounded open sets in $\mathbb R^{n_{sd}}$ with Lipschitz boundaries, where $n_{sd}$ represents the number of space dimensions. They represent the initial and the current configurations of the body, respectively. The motion of the body is described by a family of diffeomorphisms, parametrized by the time coordinate $t$,
\begin{align*}
& \bm\varphi_t(\cdot) = \bm\varphi(\cdot, t) : \Omega_{\bm X} \rightarrow \Omega_{\bm x} = \bm \varphi(\Omega_{\bm X}, t) = \bm \varphi_t(\Omega_{\bm X}), \quad \forall t \geq 0, \\
& \bm X \mapsto \bm x = \bm \varphi(\bm X, t) = \bm \varphi_t(\bm X), \quad \forall \bm X \in \Omega_{\bm X}.
\end{align*}
In the above, $\bm x$ is the current position of a material particle originally located at $\bm X$. This requires that $\bm\varphi(\bm X, 0) = \bm X$. The displacement and velocity of the material particle are defined as
\begin{align*}
\bm U := \bm \varphi(\bm X, t) - \bm \varphi(\bm X, 0) = \bm \varphi(\bm X, t) - \bm X, \quad
\bm V := \left. \frac{\partial \bm \varphi}{\partial t}\right|_{\bm X}= \left. \frac{\partial \bm U}{\partial t}\right|_{\bm X} = \frac{d\bm U}{dt}.
\end{align*}
In the definition of $\bm V$ and in what follows, $d\left( \cdot \right)/dt$ designates a total time derivative. The spatial velocity is defined as $\bm v := \bm V \circ \bm \varphi_t^{-1}$. Analogously, we define $\bm u := \bm U \circ \varphi_t^{-1}$. The deformation gradient, the Jacobian determinant, and the right Cauchy-Green tensor are defined as
\begin{align*}
\bm F := \frac{\partial \bm \varphi}{\partial \bm X}, \qquad
J := \textup{det}\left(\bm F \right), \qquad \bm C := \bm F^T \bm F.
\end{align*}
We define $\tilde{\bm F}$ and $\tilde{\bm C}$ as
\begin{align*}
\tilde{\bm F} := J^{-\frac13}\bm F, \qquad \tilde{\bm C} := J^{-\frac23}\bm C,
\end{align*}
which represent the distortional parts of $\bm F$ and $\bm C$, respectively. We denote the thermodynamic pressure of the continuum body as $p$. The mechanical behavior of an elastic material can be described by a Gibbs free energy $G(\tilde{\bm C}, p)$. In \cite{Liu2018}, it is shown that the Gibbs free energy enjoys a decoupled structure,
\begin{align*}
G(\tilde{\bm C}, p) = G_{iso}(\tilde{\bm C}) + G_{vol}(p),
\end{align*}
where $G_{iso}$ and $G_{vol}$ represent the isochoric and volumetric elastic responses. Under the isothermal condition, the energy equation is decoupled from the system, and it suffices to consider the following equations for the motion of the continuum body,
\begin{align}
\label{eq:strong_form_kinematic}
& \bm 0 = \frac{d\bm u}{dt} - \bm v, && \mbox{ in } \Omega_{\bm x}, \displaybreak[2] \\
\label{eq:strong_form_pressure}
& 0 = \beta(p) \frac{dp}{dt} + \nabla_{\bm x} \cdot \bm v &&\mbox{ in } \Omega_{\bm x},  \displaybreak[2] \\
\label{eq:strong_form_momentum}
& \bm 0 = \rho(p) \frac{d\bm v}{dt} - \nabla_{\bm x} \cdot \bm \sigma^{dev} + \nabla_{\bm x} p - \rho(p) \bm b, && \mbox{ in } \Omega_{\bm x}.
\end{align}
In the above system, the equations \eqref{eq:strong_form_kinematic} describe the kinematic relation between the displacement and the velocity, and the equations \eqref{eq:strong_form_pressure} and \eqref{eq:strong_form_momentum} describe the balance of mass and linear momentum. Let $\rho_0$ denote the density in the material configuration. The constitutive relations of the elastic material are represented in terms of the Gibbs free energy as follows,
\begin{align*}
\rho(p) :=& \left( \frac{d G_{vol}}{d p} \right)^{-1}, \quad \beta(p) := \frac{1}{\rho} \frac{d\rho}{d p} = -\frac{\partial^2 G_{vol}}{\partial p^2} / \frac{\partial G_{vol}}{\partial p},\displaybreak[2] \\
\bm \sigma^{dev} :=& J^{-1} \tilde{\bm F} \left( \mathbb P : \tilde{\bm S} \right) \tilde{\bm F}^T ,  \quad
\tilde{\bm S} := 2 \frac{\partial \left( \rho_0 G\right)}{\partial \tilde{\bm C}} = 2 \frac{\partial \left( \rho_0 G_{iso} \right)}{\partial \tilde{\bm C}}.
\end{align*}
Interested readers are referred to \cite{Liu2018} for a detailed derivation of the governing equations and the constitutive relations. The boundary $\Gamma_{\bm x} = \partial \Omega_{\bm x}$ can be partitioned into two non-overlapping subdivisions:
$
\Gamma_{\bm x} = \Gamma_{\bm x}^g \cup \Gamma_{\bm x}^h,
$
wherein $\Gamma^g_{\bm x}$ is the Dirichlet part of the boundary, and $\Gamma^h_{\bm x}$ is the Neumann part of the boundary. Boundary conditions can be stated as
\begin{align}
\label{eq:strong_form_dirichlet_u}
& \bm u = \bm g \quad \mbox{ on } \Gamma_{\bm x}^{g}, \quad
\bm v = \frac{d\bm g}{dt} \quad \mbox{ on } \Gamma_{\bm x}^{g}, \quad (\bm \sigma^{dev} - p\bm I) \bm n = \bm h \quad \mbox{ on } \Gamma_{\bm x}^{h}.
\end{align}
Given the initial data $\bm u_0$, $p_0$, and $\bm v_0$, the initial conditions can be stated as
\begin{align}
\label{eq:initial_condition_v}
& \bm u(\bm x, 0) = \bm u_0(\bm x), \quad
 p(\bm x, 0) = p_0(\bm x),  \quad
\bm v(\bm x, 0) = \bm v_0(\bm x).
\end{align}
The equations \eqref{eq:strong_form_kinematic}-\eqref{eq:initial_condition_v} constitute an initial-boundary value problem for hyper-elastodynamics.

\section{Numerical formulation}
\label{sec:numerical_formulation}
In this section, we present the numerical formulation for the strong-form problem. The spatial discretization is based on a VMS formulation \cite{Liu2018,Hughes1995}, and the temporal scheme is based on the generalized-$\alpha$ scheme \cite{Liu2018,Jansen2000}. A block factorization, originally introduced in \cite{Scovazzi2016}, is performed to consistently reduce the size of the linear algebra problem in the Newton-Raphson iterative algorithm.

\subsection{Variational multiscale formulation}
\label{subsec:spatial_discretization}
We consider a partition of $\bar{\Omega}_{\bm x}$ by $n_{el}$ non-overlapping, shape-regular elements $\Omega^e_{\bm x}$. The diameter of an element $\Omega_{\bm x}^e$ is denoted by $h^e$. The maximum diameter  of the elements is denoted as $h$, and $h \rightarrow 0$ as $n_{el} \rightarrow \infty$. Let $\mathsf P_{\mathsf k}\left( \Omega_{\bm x}^e \right)$ denote the space of complete polynomials of order $\mathsf k$ on $\bar{\Omega}^e_{\bm x}$. The finite element trial solution spaces for the displacement, pressure, and velocity are defined as
\begin{align*}
\mathcal S_{\bm u_h} &= \left\lbrace \bm u_h \mid \bm u_h(\cdot, t) \in \left( C^0(\Omega_{\bm x}) \right)^{n_{sd}}, t \in [0,T], \left. \bm u_h\right\rvert_{\Omega^e_{\bm x}} \in \left( \mathsf P_{\mathsf k}\left( \Omega_{\bm x}^e \right) \right)^{n_{sd}}, \bm u_h(\cdot,t) = \bm g \mbox{ on } \Gamma_{\bm x}^{g}  \right\rbrace , \displaybreak[2] \\
\mathcal S_{p_h} &= \left\lbrace p_h \mid p_h(\cdot, t) \in C^0(\Omega_{\bm x}), t \in [0,T], \left. p_h\right\rvert_{\Omega^e_{\bm x}} \in \mathsf P_{\mathsf k}\left( \Omega_{\bm x}^e \right)  \right\rbrace , \displaybreak[2] \\
 \mathcal S_{\bm v_h} &= \left\lbrace \bm v_h \mid \bm v_h(\cdot, t) \in \left( C^0(\Omega_{\bm x}) \right)^{n_{sd}}, t \in [0,T], \left. \bm v_h\right\rvert_{\Omega^e_{\bm x}} \in \left( \mathsf P_{\mathsf k}\left( \Omega_{\bm x}^e \right) \right)^{n_{sd}}, \bm v_h(\cdot,t) = \frac{d\bm g}{dt} \mbox{ on } \Gamma_{\bm x}^{g}  \right\rbrace ,
\end{align*}
and the corresponding test function spaces are defined as
\begin{align*}
\mathcal V_{\bm u_h} &= \left\lbrace \bm w_{\bm u_h} \mid \bm w_{\bm u_h} \in \left( C^0(\Omega_{\bm x}) \right)^{n_{sd}}, \left. \bm w_{\bm u_h}\right\rvert_{\Omega^e_{\bm x}} \in \left( \mathsf P_{\mathsf k}\left( \Omega_{\bm x}^e \right) \right)^{n_{sd}}, \bm w_{\bm u_h} = \bm 0 \mbox{ on } \Gamma_{\bm x}^{g}  \right\rbrace ,  \displaybreak[2] \\
\mathcal V_{p_h} &= \left\lbrace w_{p_h} \mid w_{p_h} \in C^0(\Omega_{\bm x}), \left. w_{p_h}\right\rvert_{\Omega^e_{\bm x}} \in \mathsf P_{\mathsf k}\left( \Omega_{\bm x}^e \right)  \right\rbrace , \displaybreak[2] \\
 \mathcal V_{\bm v_h} &= \left\lbrace \bm w_{\bm v_h} \mid \bm w_{\bm v_h} \in \left( C^0(\Omega_{\bm x}) \right)^{n_{sd}},  \left. \bm w_{\bm v_h}\right\rvert_{\Omega^e_{\bm x}} \in \left( \mathsf P_{\mathsf k}\left( \Omega_{\bm x}^e \right) \right)^{n_{sd}}, \bm w_{\bm v_h} = \bm 0 \mbox{ on } \Gamma_{\bm x}^{g}  \right\rbrace .
\end{align*}
The semi-discrete formulation can be stated as follows. Find $\bm y_h(t) := \left\lbrace  \bm u_h(t), p_h(t), \bm v_h(t)\right\rbrace^T \in \mathcal S_{\bm u_h} \times \mathcal S_{p_h} \times \mathcal S_{\bm v_h}$ such that for $t\in [0, T]$,
\begin{align}
\label{eq:kinematics_current}
& \bm 0 = \mathbf B_k\left( \bm w_{\bm u_h}; \dot{\bm y}_h, \bm y_h  \right) := \int_{\Omega_{\bm x}} \bm w_{\bm u_h} \cdot \left( \frac{d\bm u_h}{dt} - \bm v_h \right) d\Omega_{\bm x}, \displaybreak[2]\\
\label{eq:mass_current}
& 0 = \mathbf B_p\left( w_{p_h}; \dot{\bm y}_h, \bm y_h  \right) := \int_{\Omega_{\bm x}} w_{p_h} \beta(p_h) \frac{dp_h}{dt} + w_{p_h} \nabla_{\bm x} \cdot \bm v_h d\Omega_{\bm x} \nonumber \\
& \hspace{1.0cm} + \sum_{e} \int_{\Omega^e_{\bm x}} \bm \tau^e_M \nabla_{\bm x} w_{p_h} \cdot \left(\rho(p_h)\frac{d\bm v_h}{dt} - \nabla_{\bm x} \cdot \bm \sigma^{dev} + \nabla_{\bm x}p_h - \rho(p_h) \bm b \right) d\Omega_{\bm x}, \displaybreak[2] \\
\label{eq:momentum_current}
& \bm 0 = \mathbf B_m\left( \bm w_{\bm v_h}; \dot{\bm y}_h, \bm y_h  \right) := \int_{\Omega_{\bm x}} \bm w_{\bm v_h} \cdot \rho(p_h) \frac{d\bm v_h}{dt} + \nabla_{\bm x} \bm w_{\bm v_h} : \bm \sigma^{dev} - \nabla_{\bm x} \cdot \bm w_{\bm v_h} p_h \nonumber \\
& \hspace{1.0cm} - \bm w_{\bm v_h} \cdot \rho(p_h)  \bm b d\Omega_{\bm x} - \int_{\Gamma_{\bm x}^{h}} \bm w_{\bm v_h} \cdot \bm h d\Gamma_{\bm x}, 
\end{align}
for $\forall \left\lbrace  \bm w_{\bm u_h} , w_{p_h}, \bm w_{\bm v_h}\right\rbrace \in \mathcal V_{\bm u_h} \times \mathcal V_{p_h} \times \mathcal V_{\bm v_h}$, with $\dot{\bm y}_h(t) := \left\lbrace  d\bm u_h/dt, dp_h/dt, d\bm v_h/dt\right\rbrace^T$ and $\bm y_h(0) := \left\lbrace \bm u_{h0}, p_{h0}, \bm v_{h0} \right\rbrace^T$. Here $\bm u_{h0}$, $p_{h0}$, and $\bm v_{h0}$ are the $\mathcal L^2$ projections of the initial data onto the finite dimensional trial solution spaces. In the above and henceforth, the formulations for the kinematic equations, the mass equation, and the linear momentum equations are indicated by the subscripts $k$, $p$ and $m$, respectively. 

The terms involving $\bm \tau_M^e$ in \eqref{eq:mass_current} arise from the subgrid-scale modeling \cite{Liu2018}. These terms improve the stability of the Galerkin formulation without sacrificing the consistency. The design of the stabilization parameter $\bm \tau_M^e$ is the crux of the design of the VMS formulation. In this work, the following choices are made,
\begin{align*}
\bm \tau_M^e = \tau^e_M \bm I_{n_{sd}}, \quad \tau^e_M = c_m \frac{h^e}{c\rho}.
\end{align*}
In the above, $\bm I_{n_{sd}}$ is the second-order identity tensor; $c_m$ is a dimensionless parameter; $c$ is the maximum wave speed in the solid body. For compressible materials, $c$ is given by the bulk wave speed. Under the isotropic small-strain linear elastic assumption, $c = \sqrt{(\lambda + 2\mu)/\rho_0}$, where $\lambda$ and $\mu$ are the Lam\'e parameters. For incompressible materials, $c = \sqrt{\mu/\rho_0}$ is the shear wave speed. We point out that, although the choices made above are based on a simplified material model, the stabilization terms still provide an effective pressure stabilization mechanism for a range of elastic and inelastic problems \cite{Liu2018,Scovazzi2016,Zeng2017,Abboud2018,Hughes1988}. In this work, we fix $c_m$ to be $10^{-3}$ and restrict our discussion to the low-order finite element method (i.e. $\mathsf k=1$).

\subsection{Temporal discretization}
\label{subsec:temporal_discretization}
Based on the semi-discrete formulation \eqref{eq:kinematics_current}-\eqref{eq:momentum_current}, we invoke the generalized-$\alpha$ method \cite{Jansen2000} for time integration. The time interval $[0,T]$ is divided into a set of $n_{ts}$ subintervals of size $\Delta t_n := t_{n+1} - t_n$ delimited by a discrete time vector $\left\lbrace t_n \right\rbrace_{n=0}^{n_{ts}}$. The solution vector and its first-order time derivative evaluated at the time step $t_n$ are denoted as $\bm y_n$ and $\dot{\bm y}_n$; the basis function for the discrete function spaces is denoted as $N_A$. With these notations, the residual vectors can be represented as
\begin{align*}
\boldsymbol{\mathrm R}_k\left(\dot{\bm y}_{n}, \bm y_{n} \right) &: = \left\lbrace \mathbf B_k\left( N_A \bm e_i ;  \dot{\bm y}_{n}, \bm y_{n} \right) \right\rbrace , \displaybreak[2] \\
\boldsymbol{\mathrm R}_p\left(\dot{\bm y}_{n}, \bm y_{n} \right) &: =\left\lbrace \mathbf B_p\left( N_A ;  \dot{\bm y}_{n}, \bm y_{n} \right) \right\rbrace, \displaybreak[2] \\
\boldsymbol{\mathrm R}_m\left(\dot{\bm y}_{n}, \bm y_{n} \right) &: =\left\lbrace \mathbf B_m\left( N_A \bm e_i ;  \dot{\bm y}_{n}, \bm y_{n} \right) \right\rbrace.
\end{align*}
The fully discrete scheme can be stated as follows. At time step $t_n$, given $\dot{\bm y}_n$, $\bm y_n$, the time step size $\Delta t_n$, and the parameters $\alpha_m$, $\alpha_f$, and $\gamma$, find $\dot{\bm y}_{n+1}$ and $\bm y_{n+1}$ such that
\begin{align}
\label{eq:gen_alpha_fully_discrete_kinematic_eqn}
& \boldsymbol{\mathrm R}_k(\dot{\bm y}_{n+\alpha_m}, \bm y_{n+\alpha_f}) = \bm 0, \displaybreak[2] \\
\label{eq:gen_alpha_fully_discrete_pressure_rate}
& \boldsymbol{\mathrm R}_p(\dot{\bm y}_{n+\alpha_m}, \bm y_{n+\alpha_f}) = \bm 0, \displaybreak[2] \\
\label{eq:gen_alpha_fully_discrete_linear_momentum}
& \boldsymbol{\mathrm R}_m(\dot{\bm y}_{n+\alpha_m}, \bm y_{n+\alpha_f}) = \bm 0, \displaybreak[2] \\
\label{eq:gen_alpha_def_y_n_plus_1}
& \bm y_{n+1} = \bm y_{n} + \Delta t_n \dot{\bm y}_n, + \gamma \Delta t_n \left( \dot{\bm y}_{n+1} - \dot{\bm y}_{n}\right), \displaybreak[2] \\
\label{eq:gen_alpha_def_y_n_alpha_m}
& \dot{\bm y}_{n+\alpha_m} = \dot{\bm y}_{n} + \alpha_m \left(\dot{\bm y}_{n+1} - \dot{\bm y}_{n} \right), \displaybreak[2] \\
\label{eq:gen_alpha_def_y_n_alpha_f}
& \bm y_{n+\alpha_f} = \bm y_{n} + \alpha_f \left( \bm y_{n+1} - \bm y_{n} \right).
\end{align}
The choice of the parameters $\alpha_m$, $\alpha_f$ and $\gamma$ determines the accuracy and stability of the temporal scheme. Importantly, the high-frequency dissipation can be controlled via a proper parametrization of these parameters, while maintaining second-order accuracy and unconditional stability (for linear problems). For first-order dynamic problems, the parameters are chosen as
\begin{align*}
\alpha_m = \frac12 \left( \frac{3-\varrho_{\infty}}{1+\varrho_{\infty}} \right), \quad \alpha_f = \frac{1}{1+\varrho_{\infty}}, \quad \gamma = \frac{1}{1+\varrho_{\infty}},
\end{align*}
wherein $\varrho_{\infty} \in [0,1]$ denotes the spectral radius of the amplification matrix at the highest mode \cite{Jansen2000}. We adopt $\varrho_{\infty} = 0.5$ for all computations presented in this work.
\begin{remark}
Interested readers are referred to \cite{Chung1993} for the parametrization of the parameters for second-order structural dynamics. A recent study shows that using the generalized-$\alpha$ method for the first-order structural dynamics enjoys improved dissipation and dispersion properties and does not suffer from overshoot \cite{Kadapa2017}. Moreover, using a first-order structural dynamic model is quite propitious for the design of a FSI scheme \cite{Liu2018}.
\end{remark}

\subsection{A Segregated predictor multi-corrector algorithm}
\label{subsec:seg_algorithm}
One may apply an inverse of the mass matrix at both sides of the equations \eqref{eq:gen_alpha_fully_discrete_kinematic_eqn} and obtain the following simplified kinematic equations,
\begin{align}
\label{eq:gen_alpha_fully_discrete_kinematic_new}
\overline{\boldsymbol{\mathrm R}}_k(\dot{\bm y}_{n+\alpha_m}, \bm y_{n+\alpha_f}) := \dot{\bm u}_{n+\alpha_m} - \bm v_{n+\alpha_f} = \bm 0.
\end{align}
This procedure can be regarded as the application of a left preconditioner on the nonlinear algebraic equations. The new equations \eqref{eq:gen_alpha_fully_discrete_kinematic_new}, together with \eqref{eq:gen_alpha_fully_discrete_pressure_rate} and \eqref{eq:gen_alpha_fully_discrete_linear_momentum}, constitute the system of nonlinear algebraic equations to be solved in each time step. The Newton-Raphson method with consistent linearization is invoked to solve the nonlinear system of equations. At the time step $t_{n+1}$, the solution vector $\bm y_{n+1}$ is solved by means of a predictor multi-corrector algorithm. We denote $\bm y_{n+1,(l)} := \left\lbrace \bm u_{n+1,(l)}, p_{n+1,(l)}, \bm v_{n+1,(l)} \right\rbrace^T$ as the solution vector at the Newton-Raphson iteration step $l=0, \cdots, l_{max}$. The residual vectors evaluated at the iteration stage $l$ are denoted as
\begin{align*}
\boldsymbol{\mathrm R}_{(l)} &:= \left\lbrace \overline{\boldsymbol{\mathrm R}}_{k,(l)}, \boldsymbol{\mathrm R}_{p,(l)}, \boldsymbol{\mathrm R}_{m,(l)} \right\rbrace^T, \\
\overline{\boldsymbol{\mathrm R}}_{k,(l)} &:= \overline{\boldsymbol{\mathrm R}}_k\left( \dot{\bm y}_{n+\alpha_m, (l)}, \bm y_{n+\alpha_f, (l)} \right), \nonumber \\
\boldsymbol{\mathrm R}_{p,(l)} &:= \boldsymbol{\mathrm R}_p\left( \dot{\bm y}_{n+\alpha_m, (l)}, \bm y_{n+\alpha_f, (l)} \right), \nonumber \\
\boldsymbol{\mathrm R}_{m,(l)} &:= \boldsymbol{\mathrm R}_m\left( \dot{\bm y}_{n+\alpha_m, (l)}, \bm y_{n+\alpha_f, (l)} \right). \nonumber
\end{align*}
The consistent tangent matrix associated with the above residual vectors is
\begin{align*}
\boldsymbol{\mathrm K}_{(l)} =
\begin{bmatrix}
\boldsymbol{\mathrm K}_{k,(l), \dot{\bm u}} & \boldsymbol{\mathrm K}_{k,(l), \dot{p}} & \boldsymbol{\mathrm K}_{k,(l), \dot{\bm v}} \\[0.3mm]
\boldsymbol{\mathrm K}_{p,(l), \dot{\bm u}} & \boldsymbol{\mathrm K}_{p,(l), \dot{p}} & \boldsymbol{\mathrm K}_{p,(l), \dot{\bm v}} \\[0.3mm]
\boldsymbol{\mathrm K}_{m,(l), \dot{\bm u}} & \boldsymbol{\mathrm K}_{m,(l), \dot{p}} & \boldsymbol{\mathrm K}_{m,(l), \dot{\bm v}}
\end{bmatrix},
\end{align*}
wherein
\begin{align*}
& \boldsymbol{\mathrm K}_{k,(l), \dot{\bm u}} := \alpha_m \frac{\partial \overline{\boldsymbol{\mathrm R}}_{k,(l)}\left( \dot{\bm y}_{n+\alpha_m, (l)}, \bm y_{n+\alpha_f, (l)} \right)}{\partial \dot{\bm u}_{n+\alpha_m}} = \alpha_m \boldsymbol{\mathrm I}, \displaybreak[2] \\
& \boldsymbol{\mathrm K}_{k,(l), \dot{p}} := \bm 0, \displaybreak[2] \\
& \boldsymbol{\mathrm K}_{k,(l), \dot{\bm v}} := \alpha_f \gamma \Delta t_n \frac{\partial \overline{\boldsymbol{\mathrm R}}_{k,(l)}\left( \dot{\bm y}_{n+\alpha_m, (l)}, \bm y_{n+\alpha_f, (l)} \right)}{\partial \bm v_{n+\alpha_f}} = -\alpha_f \gamma \Delta t_n \boldsymbol{\mathrm I}.
\end{align*}
As was realized in \cite{Scovazzi2016}, this special block structure in the first row of $\boldsymbol{\mathrm K}_{(l)}$ can be utilized for a block factorization,
\begin{align}
\label{eq:K_block_decomposition}
\boldsymbol{\mathrm K}_{(l)} =& \begin{bmatrix}
\boldsymbol{\mathrm K}_{k,(l), \dot{\bm u}} & \boldsymbol{\mathrm K}_{k,(l), \dot{p}} & \boldsymbol{\mathrm K}_{k,(l), \dot{\bm v}} \\[0.3mm]
\boldsymbol{\mathrm K}_{p,(l), \dot{\bm u}} & \boldsymbol{\mathrm K}_{p,(l), \dot{p}} & \boldsymbol{\mathrm K}_{p,(l), \dot{\bm v}} \\[0.3mm]
\boldsymbol{\mathrm K}_{m,(l), \dot{\bm u}} & \boldsymbol{\mathrm K}_{m,(l), \dot{p}} & \boldsymbol{\mathrm K}_{m,(l), \dot{\bm v}} 
\end{bmatrix} =
\begin{bmatrix}
\alpha_m \boldsymbol{\mathrm I} & \bm 0 & -\alpha_f \gamma \Delta t_n \boldsymbol{\mathrm I} \\[0.3mm]
\boldsymbol{\mathrm K}_{p,(l), \dot{\bm u}} & \boldsymbol{\mathrm K}_{p,(l), \dot{p}} & \boldsymbol{\mathrm K}_{p,(l), \dot{\bm v}} \\[0.3mm]
\boldsymbol{\mathrm K}_{m,(l), \dot{\bm u}} & \boldsymbol{\mathrm K}_{m,(l), \dot{p}} & \boldsymbol{\mathrm K}_{m,(l), \dot{\bm v}} 
\end{bmatrix}
\nonumber \\
=&
\begin{bmatrix}
\boldsymbol{\mathrm I} & \bm 0 & \bm 0 \\[0.3mm]
\frac{1}{\alpha_m}\boldsymbol{\mathrm K}_{p,(l), \dot{\bm u}} & \boldsymbol{\mathrm K}_{p,(l), \dot{p}} & \boldsymbol{\mathrm K}_{p,(l), \dot{\bm v}} + \frac{\alpha_f \gamma \Delta t_n}{\alpha_m} \boldsymbol{\mathrm K}_{p,(l), \dot{\bm u}} \\[0.3mm]
\frac{1}{\alpha_m}\boldsymbol{\mathrm K}_{m,(l), \dot{\bm u}} & \boldsymbol{\mathrm K}_{m,(l), \dot{p}} & \boldsymbol{\mathrm K}_{m,(l), \dot{\bm v}} +  \frac{\alpha_f \gamma \Delta t_n}{\alpha_m} \boldsymbol{\mathrm K}_{m,(l), \dot{\bm u}}
\end{bmatrix}
\begin{bmatrix}
\alpha_m \boldsymbol{\mathrm I} & \bm 0 & -\alpha_f \gamma \Delta t_n \boldsymbol{\mathrm I} \\[0.3mm]
\bm 0 & \boldsymbol{\mathrm I} & \bm 0 \\[0.3mm]
\bm 0 & \bm 0 & \boldsymbol{\mathrm I}
\end{bmatrix}.
\end{align}
With \eqref{eq:K_block_decomposition}, the solution procedure of the linear system of equations in the Newton-Raphson method can be consistently reduced to a two-stage algorithm \cite{Liu2018,Scovazzi2016,Rossi2016}. In the first stage, one obtains the increments of the pressure and velocity at the iteration step $l$ by solving the following linear system,
\begin{align}
\label{eq:NR_smaller_eqn_1}
&
\begin{bmatrix}
\boldsymbol{\mathrm K}_{m,(l), \dot{\bm v}} +  \frac{\alpha_f \gamma \Delta t_n}{\alpha_m} \boldsymbol{\mathrm K}_{m,(l), \dot{\bm u}} & \boldsymbol{\mathrm K}_{m,(l), \dot{p}} \\[0.3mm]
\boldsymbol{\mathrm K}_{p,(l), \dot{\bm v}} + \frac{\alpha_f \gamma \Delta t_n}{\alpha_m} \boldsymbol{\mathrm K}_{p,(l), \dot{\bm u}} & \boldsymbol{\mathrm K}_{p,(l), \dot{p}}
\end{bmatrix}
\begin{bmatrix}
\Delta \dot{\bm v}_{n+1,(l)} \\[0.3mm]
\Delta \dot{p}_{n+1,(l)}
\end{bmatrix}
= -
\begin{bmatrix}
\boldsymbol{\mathrm R}_{m,(l)}- \frac{1}{\alpha_m}  \boldsymbol{\mathrm K}_{m,(l), \dot{\bm u}}\overline{\boldsymbol{\mathrm R}}_{k,(l)} \\[0.3mm]
\boldsymbol{\mathrm R}_{p,(l)} - \frac{1}{\alpha_m}  \boldsymbol{\mathrm K}_{p,(l), \dot{\bm u}}\overline{\boldsymbol{\mathrm R}}_{k,(l)}
\end{bmatrix}.
\end{align}
In the second stage, one obtains the increments for the displacement by
\begin{align}
\label{eq:NR_smaller_eqn_2}
& \Delta \dot{\bm u}_{n+1,(l)} = \frac{\alpha_f \gamma \Delta t_n}{\alpha_m} \Delta \dot{\bm v}_{n+1,(l)} - \frac{1}{\alpha_m}\overline{\boldsymbol{\mathrm R}}^k_{(l)}.
\end{align}
To simplify notations in the following discussion, we denote 
\begin{align}
\label{eq:def_matrix_A}
\boldsymbol{\mathrm A}_{(l)} :=& \boldsymbol{\mathrm K}_{m,(l), \dot{\bm v}} +  \frac{\alpha_f \gamma \Delta t_n}{\alpha_m} \boldsymbol{\mathrm K}_{m,(l), \dot{\bm u}}, \quad
\boldsymbol{\mathrm B}_{(l)} := \boldsymbol{\mathrm K}_{m,(l), \dot{p}}, \\
\label{eq:def_matrix_D}
\boldsymbol{\mathrm C}_{(l)} :=& \boldsymbol{\mathrm K}_{p,(l), \dot{\bm v}} + \frac{\alpha_f \gamma \Delta t_n}{\alpha_m} \boldsymbol{\mathrm K}_{p,(l), \dot{\bm u}}, \quad
\boldsymbol{\mathrm D}_{(l)} := \boldsymbol{\mathrm K}_{p,(l), \dot{p}}.
\end{align}
\begin{remark}
In \cite{Liu2018}, it was shown that $\overline{\boldsymbol{\mathrm R}}^k_{(l)} = \bm 0$ for $l \geq 2$ for general predictor multi-corrector algorithms; in \cite{Rossi2016}, a special predictor is chosen so that $\overline{\boldsymbol{\mathrm R}}^k_{(l)} = \bm 0$ for $l \geq 1$. 
\end{remark}
\begin{remark}
\label{remark:structure_of_A_B_C}
In \ref{app:consistent_linearization}, the detailed formula for the block matrices are given, and it can be observed that $\boldsymbol{\mathrm A}_{(l)}$ consists primarily of a mass matrix and a stiffness matrix; $\boldsymbol{\mathrm B}_{(l)}$ is a discrete gradient operator; $\boldsymbol{\mathrm C}_{(l)}$ is dominated by a discrete divergence operator; $\boldsymbol{\mathrm D}_{(l)}$ contains a mass matrix scaled with $\beta$ and contributions from the stabilization terms.
\end{remark}
Based on the above discussion, a predictor multi-corrector algorithm for solving the nonlinear algebraic equations in each time step can be summarized as follows.

\noindent \textbf{Predictor stage}: Set:
\begin{align*}
\bm y_{n+1, (0)} = \bm y_{n}, \quad
\dot{\bm y}_{n+1, (0)} = \frac{\gamma - 1}{\gamma} \dot{\bm y}_{n}.
\end{align*}

\noindent \textbf{Multi-corrector stage}:
Repeat the following steps for \(l = 1, \dots, l_{max}\):
\begin{enumerate}
\item Evaluate the solution vectors at the intermediate stages:
\begin{align*}
\dot{\bm y}_{n+\alpha_m, (l)} &= \dot{\bm y}_n + \alpha_m \left( \dot{\bm y}_{n+1,(l-1)} - \dot{\bm y}_n \right), \\
\bm y_{n+\alpha_f, (l)} &= \bm y_n + \alpha_f \left( \bm y_{n+1,(l-1)} - \bm y_n \right).
\end{align*}

\item Assemble the residual vectors $\boldsymbol{\mathrm R}_{m,(l)}$ and $\boldsymbol{\mathrm R}_{p,(l)}$ using the solution evaluated at the intermediate stages.

\item Let $\|\boldsymbol{\mathrm R}_{(l)}\|_{\mathfrak{l}^2}$ denote the $\mathfrak l^2$-norm of the residual vector. If either one of the following stopping criteria 
\begin{align*}
& \frac{\|\boldsymbol{\mathrm R}_{(l)}\|_{\mathfrak l^2}}{\|\boldsymbol{\mathrm R}_{(0)}\|_{\mathfrak l^2}} \leq \textup{tol}_{\textup{R}}, \qquad \|\boldsymbol{\mathrm R}_{(l)}\|_{\mathfrak l^2} \leq \textup{tol}_{\textup{A}},
\end{align*}
is satisfied for two prescribed tolerances $\textup{tol}_{\textup{R}}$, $\textup{tol}_{\textup{A}}$, set the solution vector at time step $t_{n+1}$ as $\dot{\bm y}_{n+1} = \dot{\bm y}_{n+1, (l-1)}$ and $\bm y_{n+1} = \bm y_{n+1, (l-1)}$, and exit the multi-corrector stage; otherwise, continue to step 4.

\item Assemble the tangent matrices \eqref{eq:def_matrix_A}-\eqref{eq:def_matrix_D}.

\item Solve the following linear system of equations for $\Delta \dot{p}_{n+1,(l)}$ and $\Delta \dot{\bm v}_{n+1,(l)}$,
\begin{align}
\label{eq:pre-multi-corrector-stage-5-matrix-problems}
\begin{bmatrix}
\boldsymbol{\mathrm A}_{(l)} & \boldsymbol{\mathrm B}_{(l)} \\[0.3mm]
\boldsymbol{\mathrm C}_{(l)} & \boldsymbol{\mathrm D}_{(l)}
\end{bmatrix}
\begin{bmatrix}
\Delta \dot{\bm v}_{n+1,(l)} \\[0.3mm]
\Delta \dot{p}_{n+1,(l)}
\end{bmatrix}
= -
\begin{bmatrix}
\boldsymbol{\mathrm R}_{m,(l)} \\[0.3mm]
\boldsymbol{\mathrm R}_{p,(l)}
\end{bmatrix}.
\end{align}

\item Obtain $\Delta \dot{\bm u}_{n+1,(l)}$ from the relation \eqref{eq:NR_smaller_eqn_2}.

\item Update the solution vector as
\begin{align*}
\dot{\bm y}_{n+1,(l)} &= \dot{\bm y}_{n+1,(l)} + \Delta \dot{\bm y}_{n+1,(l)}, \\
\bm y_{n+1,(l)} &= \bm y_{n+1,(l)} + \gamma \Delta t_n \Delta \dot{\bm y}_{n+1,(l)}.
\end{align*}
and return to step 1.
\end{enumerate}
For all the numerical simulations presented in this work, we adopt the tolerances for the nonlinear iteration as $\textup{tol}_{\textup{R}} = \textup{tol}_{\textup{A}} = 10^{-6}$ and the maximum number of iterations as $l_{max}=20$.

\section{Iterative linear solver}
\label{sec:linear_solver}
In the predictor multi-corrector algorithm presented above, the linear system of equations \eqref{eq:pre-multi-corrector-stage-5-matrix-problems} is solved repeatedly, and this step constitutes the major cost for implicit dynamic calculations. In this section, we design an iterative solution procedure for the linear problem $\mathcal A \bm x = \bm r$, in which the matrix and vectors adopt the following block structure,
\begin{align*}
\mathcal A :=
\begin{bmatrix}
\boldsymbol{\mathrm A} & \boldsymbol{\mathrm B} \\[0.3mm]
\boldsymbol{\mathrm C} & \boldsymbol{\mathrm D}
\end{bmatrix}, \quad
\bm x :=
\begin{bmatrix}
\bm x_{\bm v} \\[0.3em]
\bm x_{p}
\end{bmatrix},
\quad
\bm r := 
\begin{bmatrix}
\bm r_{\bm v} \\[0.3em]
\bm r_{p}
\end{bmatrix}.
\end{align*}

Since its inception, GMRES is among the most popular iterative methods for solving sparse nonsymmetric matrix problems. With a proper preconditioner $\mathcal P$, the convergence rate of iterative methods like GMRES can be significantly expedited. Roughly speaking, in the GMRES iteration, one constructs the Krylov subspace and search for the solution that minimize the residual in this Krylov subspace by the Arnoldi algorithm \cite{Saad1986,Saad1993}. To construct the Krylov subspace, one applies $\mathcal A \mathcal P^{-1}$ to the residual vector in order to enlarge the Krylov subspace. This procedure corresponds to first solving a linear system of equations associated with $\mathcal P$ and then performing a matrix-vector multiplication associated with $\mathcal A$. Often times, to reduce the computational burden, the GMRES algorithm is restarted every $\mathfrak m$ steps. Within this work, this algorithm is denoted as GMRES($\mathfrak m$).

In Section \ref{subsec:sym_diag_scale}, we perform a diagonal scaling for $\mathcal A$ with the purpose of improving the condition number \cite{May2008,Moghadam2015,Shakib1989}. In Section \ref{subsec:block_pc}, we introduce the block factorization of $\mathcal A$ and present the SCR algorithm. In Section \ref{subsec:coupled_approach}, we present the coupled approach with a particular focus on the SIMPLE preconditioner. In Section \ref{subsec:FGMRES_nested_pc}, the nested block preconditioning technique is introduced as a combination of the SCR approach and the coupled approach.

\subsection{Symmetrically diagonal scaling}
\label{subsec:sym_diag_scale}
Before constructing an iterative method, we first apply a symmetrically diagonal scaling to the matrix $\mathcal A$. This approach is adopted to improve the condition number of the matrix problem and is sometimes referred to as a ``pre-preconditioning" technique \cite{Shakib1989}. We introduce $\mathcal W$ as a diagonal matrix defined as follows,
\begin{align*}
\mathcal W_{ii} := 
\begin{cases}
\left( |\mathcal A_{ii}| \right)^{-\frac12},  & \quad \text{if } |\mathcal A_{ii}| \geq \epsilon_{diag} \\
1.0,  & \quad \text{if } |\mathcal A_{ii}| < \epsilon_{diag}
\end{cases}
.
\end{align*}
In the above definition, $\epsilon_{diag}$ is a user-specified tolerance to avoid undefined or unstable numerical operations. In this work, we set $\epsilon_{diag} = 1.0\times 10^{-15}$. Applying $\mathcal W$ as a left and right preconditioner simultaneously, we obtain an altered system as
\begin{align}
\label{eq:sym_diag_scaled_linear_system}
\mathcal A^* \bm x^* = \bm r^*,
\end{align}
wherein $\mathcal A^* := \mathcal W \mathcal A \mathcal W$, $\bm x^* := \mathcal W^{-1}\bm x$, and $\bm r^* = \mathcal W \bm r$. The iterative methods discussed in the subsequent sections are applied to the above system. Once $x^*$ is obtained from \eqref{eq:sym_diag_scaled_linear_system}, one has to perform $\bm x = \mathcal W \bm x^*$ to recover the true solution. In the remainder of Section \ref{sec:linear_solver}, we focus on solving \eqref{eq:sym_diag_scaled_linear_system}, and for notational simplicity, the superscript $*$ is neglected.

\subsection{Schur complement reduction}
\label{subsec:block_pc}
Recall that $\mathcal A$ adopts the block factorization
\begin{align}
\label{eq:a_ldu_block_factorization}
\mathcal A = \mathcal L \mathcal D \mathcal U =
\begin{bmatrix}
\boldsymbol{\mathrm I} & \boldsymbol{\mathrm O} \\[0.3em]
\boldsymbol{\mathrm C} \boldsymbol{\mathrm A}^{-1} & \boldsymbol{\mathrm I}
\end{bmatrix}
\begin{bmatrix}
\boldsymbol{\mathrm A} & \boldsymbol{\mathrm O} \\[0.3em]
\boldsymbol{\mathrm O} & \boldsymbol{\mathrm S}
\end{bmatrix}
\begin{bmatrix}
\boldsymbol{\mathrm I} & \boldsymbol{\mathrm A}^{-1} \boldsymbol{\mathrm B} \\[0.3em]
\boldsymbol{\mathrm O} & \boldsymbol{\mathrm I}
\end{bmatrix},
\end{align}
wherein $\boldsymbol{\mathrm I}$ is the identity matrix, $\boldsymbol{\mathrm O}$ is the zero matrix, and $\boldsymbol{\mathrm S} := \boldsymbol{\mathrm D} - \boldsymbol{\mathrm C} \boldsymbol{\mathrm A}^{-1} \boldsymbol{\mathrm B}$ is the Schur complement of $\boldsymbol{\mathrm A}$. Applying $\mathcal L^{-1}$ on both sides of the equation $\mathcal A \bm x = \bm r$, one obtains
\begin{align*}
\begin{bmatrix}
\boldsymbol{\mathrm A} & \boldsymbol{\mathrm B} \\[0.3em]
\boldsymbol{\mathrm O} & \boldsymbol{\mathrm S}
\end{bmatrix}
\begin{bmatrix}
\bm x_{\bm v} \\[0.3em]
\bm x_{p}
\end{bmatrix}
&=
\begin{bmatrix}
\boldsymbol{\mathrm I} & \boldsymbol{\mathrm O} \\[0.3em]
\boldsymbol{\mathrm C} \boldsymbol{\mathrm A}^{-1} & \boldsymbol{\mathrm I}
\end{bmatrix}^{-1}
\begin{bmatrix}
\bm r_{\bm v} \\[0.3em]
\bm r_{p}
\end{bmatrix}
=
\begin{bmatrix}
\boldsymbol{\mathrm I} & \boldsymbol{\mathrm O} \\[0.3em]
-\boldsymbol{\mathrm C} \boldsymbol{\mathrm A}^{-1} & \boldsymbol{\mathrm I}
\end{bmatrix}
\begin{bmatrix}
\bm r_{\bm v} \\[0.3em]
\bm r_{p}
\end{bmatrix} =
\begin{bmatrix}
\bm r_{\bm v} \\[0.3em]
\bm r_{p} - \boldsymbol{\mathrm C} \boldsymbol{\mathrm A}^{-1} \bm r_{\bm v}
\end{bmatrix}.
\end{align*}
The upper triangular block matrix problem can be solved by a back substitution. Consequently, the solution procedure for $\mathcal A \bm x = \bm r$ can be summarized as the following segregated algorithm \cite{May2008,Lun2010,Furuichi2011}.
\begin{algorithm}[H]
\caption{Solution procedure for $\mathcal A \bm x = \bm r$ based on SCR.}
\label{algorithm:exact_block_factorization}
\begin{algorithmic}[1]
\\ Solve for an intermediate velocity $\hat{\bm x}_{\bm v}$ from the equation
\begin{align}
\label{eq:seg_sol_int_disp}
\boldsymbol{\mathrm A} \hat{\bm x}_{\bm v} = \bm r_{\bm v}.
\end{align}
\\ Update the continuity residual by $\bm r_{p} \gets  \bm r_{p} - \boldsymbol{\mathrm C} \hat{\bm x}_{\bm v}$.
\\ Solve for $\bm x_p$ from the equation
\begin{align}
\label{eq:seg_sol_pres}
\boldsymbol{\mathrm S} \bm x_p = \bm r_{p}.
\end{align}
\\ Update the momentum residual by $\bm r_{\bm v} \gets \bm r_{\bm v} - \boldsymbol{\mathrm B} \bm x_{p}$.
\\ Solve for $\bm x_{\bm v}$ from the equation
\begin{align}
\label{eq:seg_sol_disp}
\boldsymbol{\mathrm A} \bm x_{\bm v} = \bm r_{\bm v}.
\end{align}
\end{algorithmic}
\end{algorithm}
For hyper-elastodynamics problems, it is reasonable to apply GMRES preconditioned by AMG for \eqref{eq:seg_sol_int_disp} and \eqref{eq:seg_sol_disp}. The stopping condition for solving with $\boldsymbol{\mathrm A}$ includes the tolerance for the relative error $\delta^r_A$, the tolerance for the absolute error $\delta^a_A$, and the maximum number of iterations $n^{max}_A$. In \eqref{eq:seg_sol_pres}, the Schur complement is a dense matrix due to the presence of $\boldsymbol{\mathrm A}^{-1}$ in its definition. It is expensive and often impossible to directly compute with $\boldsymbol{\mathrm S}$. Recall that in a Krylov subspace method, the search space is iteratively expanded by performing matrix-vector multiplications. Although the algebraic form of $\boldsymbol{\mathrm S}$ is impractical to obtain, its action on a vector is readily available through the following ``matrix-free" algorithm \cite{May2008,Furuichi2011}.
\begin{algorithm}[H]
\caption{The multiplication of $\boldsymbol{\mathrm S}$ with a vector $\bm x_p$.}
\label{algorithm:matrix_free_mat_vec_for_S}
\begin{algorithmic}[1]
\State Compute the matrix-vector multiplication $\hat{\bm x}_p \gets \boldsymbol{\mathrm D} \bm x_p$.
\State Compute the matrix-vector multiplication $\bar{\bm x}_p \gets \boldsymbol{\mathrm B} \bm x_p$.
\State Solve for $\tilde{\bm x}_p$ from the linear system 
\begin{align}
\label{eq:S_inner_A_eqn}
\boldsymbol{\mathrm A} \tilde{\bm x}_p = \bar{\bm x}_p.
\end{align}
\State Compute the matrix-vector multiplication $\bar{\bm x}_p \gets \boldsymbol{\mathrm C} \tilde{\bm x}_p$.
\State \Return $\hat{\bm x}_p - \bar{\bm x}_p$.
\end{algorithmic}
\end{algorithm}
In Algorithm \ref{algorithm:matrix_free_mat_vec_for_S}, the action of $\boldsymbol{\mathrm S}$ on a vector is realized through a series of matrix-vector multiplications, and the action of $\boldsymbol{\mathrm A}^{-1}$ on a vector is achieved by solving the linear system \eqref{eq:S_inner_A_eqn}. This solver is located inside the solution procedure of \eqref{eq:seg_sol_pres}, and we call it the \textit{inner solver}. The stopping condition of the inner solver includes the tolerance for the relative error $\delta^r_I$, the tolerance for the absolute error $\delta^a_I$, and the maximum number of iterations $n^{max}_I$.

With Algorithm \ref{algorithm:matrix_free_mat_vec_for_S}, one can construct a Krylov subspace for $\boldsymbol{\mathrm S}$ and solve the equation \eqref{eq:seg_sol_pres}. However, without preconditioning, GMRES may stagnate or even break down. More importantly, each matrix-vector multiplication given in Algorithm \ref{algorithm:matrix_free_mat_vec_for_S} involves solving a linear system \eqref{eq:S_inner_A_eqn}, and this inevitably makes the matrix-vector multiplication quite expensive. To mitigate the number of this expensive matrix-vector multiplications, we solve \eqref{eq:seg_sol_pres} with $\hat{\boldsymbol{\mathrm S}}:= \boldsymbol{\mathrm D} - \boldsymbol{\mathrm C} \left(\textup{diag}\left(\boldsymbol{\mathrm A}\right)\right)^{-1} \boldsymbol{\mathrm B}$ as a right preconditioner \cite{Berger-Vergiat2016}. If the time step size is small, $\boldsymbol{\mathrm A}$ is dominated by the mass matrix, and $\hat{\boldsymbol{\mathrm S}}$ acts as an effective preconditioner for solving \eqref{eq:seg_sol_pres}. On the other side, if the time step size is large, $\boldsymbol{\mathrm A}$ is dominated by the stiffness matrix. The situation then is analogous to the Stokes problem, where the Schur complement is spectrally equivalent to an identity matrix. We may reasonably expect that an unpreconditioned GMRES using Algorithm \ref{algorithm:matrix_free_mat_vec_for_S} is sufficient for solving \eqref{eq:seg_sol_pres}. Still, using $\hat{\boldsymbol{\mathrm S}}$ may accelerate the convergence rate. Therefore, we solve \eqref{eq:seg_sol_pres} by GMRES, where the stopping criteria include the tolerance for the relative error $\delta^r_S$, the tolerance for the absolute error $\delta^a_S$, and the maximum number of iterations $n^{max}_S$.


\subsection{Coupled approach with block preconditioners}
\label{subsec:coupled_approach}
The block factorization \eqref{eq:a_ldu_block_factorization} also inspires the design of a preconditioner for $\mathcal A$. Following the nomenclature used in \cite{Quarteroni2000}, we use $\boldsymbol{\mathrm H}_1$ and $\boldsymbol{\mathrm H}_2$ to denote the approximations of $\boldsymbol{\mathrm A}^{-1}$ in the Schur complement and the upper triangular matrix $\mathcal U$, respectively. This results in a block preconditioner expressed as
\begin{align*}
\hat{\mathcal P} = 
\begin{bmatrix}
\boldsymbol{\mathrm I} & \boldsymbol{\mathrm O} \\[0.3em]
\boldsymbol{\mathrm C} \boldsymbol{\mathrm A}^{-1} & \boldsymbol{\mathrm I}
\end{bmatrix}
\begin{bmatrix}
\boldsymbol{\mathrm A} & \boldsymbol{\mathrm O} \\[0.3em]
\boldsymbol{\mathrm O} & \boldsymbol{\mathrm D} - \boldsymbol{\mathrm C} \boldsymbol{\mathrm H}_1 \boldsymbol{\mathrm B}
\end{bmatrix}
\begin{bmatrix}
\boldsymbol{\mathrm I} & \boldsymbol{\mathrm H}_2 \boldsymbol{\mathrm B} \\[0.3em]
\boldsymbol{\mathrm O} & \boldsymbol{\mathrm I}
\end{bmatrix} =
\begin{bmatrix}
\boldsymbol{\mathrm A} & \boldsymbol{\mathrm A}\boldsymbol{\mathrm H}_2\boldsymbol{\mathrm B} \\[0.3em]
\boldsymbol{\mathrm C} & \boldsymbol{\mathrm D} - \boldsymbol{\mathrm C} \left( \boldsymbol{\mathrm H}_1 - \boldsymbol{\mathrm H}_2\right) \boldsymbol{\mathrm B}
\end{bmatrix}.
\end{align*}
The two approximated sparse matrices are introduced so that the spectrum of $\mathcal A \hat{\mathcal P}^{-1}$ has a clustering around $\{1\}$. With the block preconditioner, one can apply the Krylov subspace method directly to solve $\mathcal A \bm x = \bm r$, and the bases of the Krylov subspace are constructed by applying $\mathcal A \hat{\mathcal P}^{-1}$ on a vector. The action of $\hat{\mathcal P}^{-1}$ is achieved through a procedure similar to the Algorithm \ref{algorithm:exact_block_factorization}. The differences are that the inner solver is not needed and one does not need to solve the equations associated with the sub-matrices to a high precision. The Krylov subspace method is typically used with a multigrid \cite{Elman2008,Cyr2012} or a domain decomposition \cite{Deparis2016} preconditioner to solve with the sub-matrices. Consequently, the algebraic definition of $\hat{\mathcal P}$ varies over iterations, and one has to apply a flexible method, like the Flexible GMRES (FGMRES) \cite{Saad1993}, as the iterative method for $\mathcal A$. Choosing $\boldsymbol{\mathrm H}_1 = \boldsymbol{\mathrm H}_2 = \textup{diag}\left(\boldsymbol{\mathrm A}\right)^{-1}$ leads to the SIMPLE preconditioner $\hat{\mathcal P}_{\textup{SIMPLE}}$ \cite{Elman2008,Quarteroni2000},
\begin{align*}
\hat{\mathcal P}_{\textup{SIMPLE}} := 
\begin{bmatrix}
\boldsymbol{\mathrm I} & \boldsymbol{\mathrm O} \\[0.3em]
\boldsymbol{\mathrm C} \boldsymbol{\mathrm A}^{-1} & \boldsymbol{\mathrm I}
\end{bmatrix}
\begin{bmatrix}
\boldsymbol{\mathrm A} & \boldsymbol{\mathrm O} \\[0.3em]
\boldsymbol{\mathrm O} & \hat{ \boldsymbol{\mathrm S} }
\end{bmatrix}
\begin{bmatrix}
\boldsymbol{\mathrm I} & \left(\textup{diag}\left(\boldsymbol{\mathrm A}\right)\right)^{-1} \boldsymbol{\mathrm B} \\[0.3em]
\boldsymbol{\mathrm O} & \boldsymbol{\mathrm I}
\end{bmatrix}
=
\begin{bmatrix}
\boldsymbol{\mathrm A} & \boldsymbol{\mathrm A}\textup{diag}\left(\boldsymbol{\mathrm A}\right)^{-1}\boldsymbol{\mathrm B} \\[0.3em]
\boldsymbol{\mathrm C} & \boldsymbol{\mathrm D}
\end{bmatrix}.
\end{align*}
The SIMPLE preconditioner is an algebraic analogue of the Semi-Implicit Method for Pressure Linked Equations (SIMPLE) \cite{Patankar1983}. It introduces a perturbation to the pressure operator in the linear momentum equation. This preconditioner and its variants are among the most popular choices for problems in CFD \cite{Cyr2012,Deparis2014}, FSI \cite{Deparis2016}, and multiphysics problems \cite{Verdugo2016,White2011}.

\begin{remark}
There are cases when the symmetry of $\boldsymbol{\mathrm A}$ is broken, and using the SIMPLE-type preconditioner leads to poor performance. It is the case in CFD with large Reynolds numbers. To take into account of the off-diagonal entries of $\boldsymbol{\mathrm A}$, sophisticated preconditioners, like the LSC preconditioner \cite{Elman2006}, have been developed. Those preconditioners have been shown to be robust with respect to the Reynolds number using inf-sup stable discretizations of the CFD problem (i.e., $\boldsymbol{\mathrm D} = \boldsymbol{\mathrm O}$). Note that, for the stabilized methods, the LSC preconditioner may not converge \cite{Cyr2012}.
\end{remark}

\subsection{Flexible GMRES algorithm with a nested block preconditioner}
\label{subsec:FGMRES_nested_pc}
The SIMPLE preconditioner can be viewed as the SCR approach built based on an inexact block factorization. Its main advantage is that the application of this preconditioner is inexpensive. However, for certain problems, this inexact factorization misses some key information of the original matrix, and stagnation of the solver is observed. We want to leverage the robustness of the SCR approach built from the exact block factorization by using it as a right preconditioner, denoted as $\hat{\mathcal P}_{SCR}$. The action of $\hat{\mathcal P}_{SCR}^{-1}$ on a vector is given by Algorithm \ref{algorithm:exact_block_factorization}, in which the equations \eqref{eq:seg_sol_int_disp}-\eqref{eq:seg_sol_disp} are solved with prescribed tolerances. The algebraic form of $\hat{\mathcal P}_{SCR}$ is defined implicitly through the solvers in Algorithm \ref{algorithm:exact_block_factorization} and varies over iterations. Assuming that the three equations \eqref{eq:seg_sol_int_disp}-\eqref{eq:seg_sol_disp} are solved exactly, the spectrum of $\mathcal A \hat{\mathcal P}_{SCR}^{-1}$ will be $\{1\}$, and the solver will converge in one iteration. Because the preconditioner varies over iterations, we invoke FGMRES as the iterative method for $\mathcal A \bm x = \bm r$. The stopping condition of the FGMRES algorithm includes the tolerance for the absolute error $\delta^a$, the tolerance for the relative error $\delta^r$, and the maximum number of iterations $n^{max}$. 

\begin{figure}
\begin{center}
\includegraphics[angle=0, trim=150 60 70 20, clip=true, scale = 0.4]{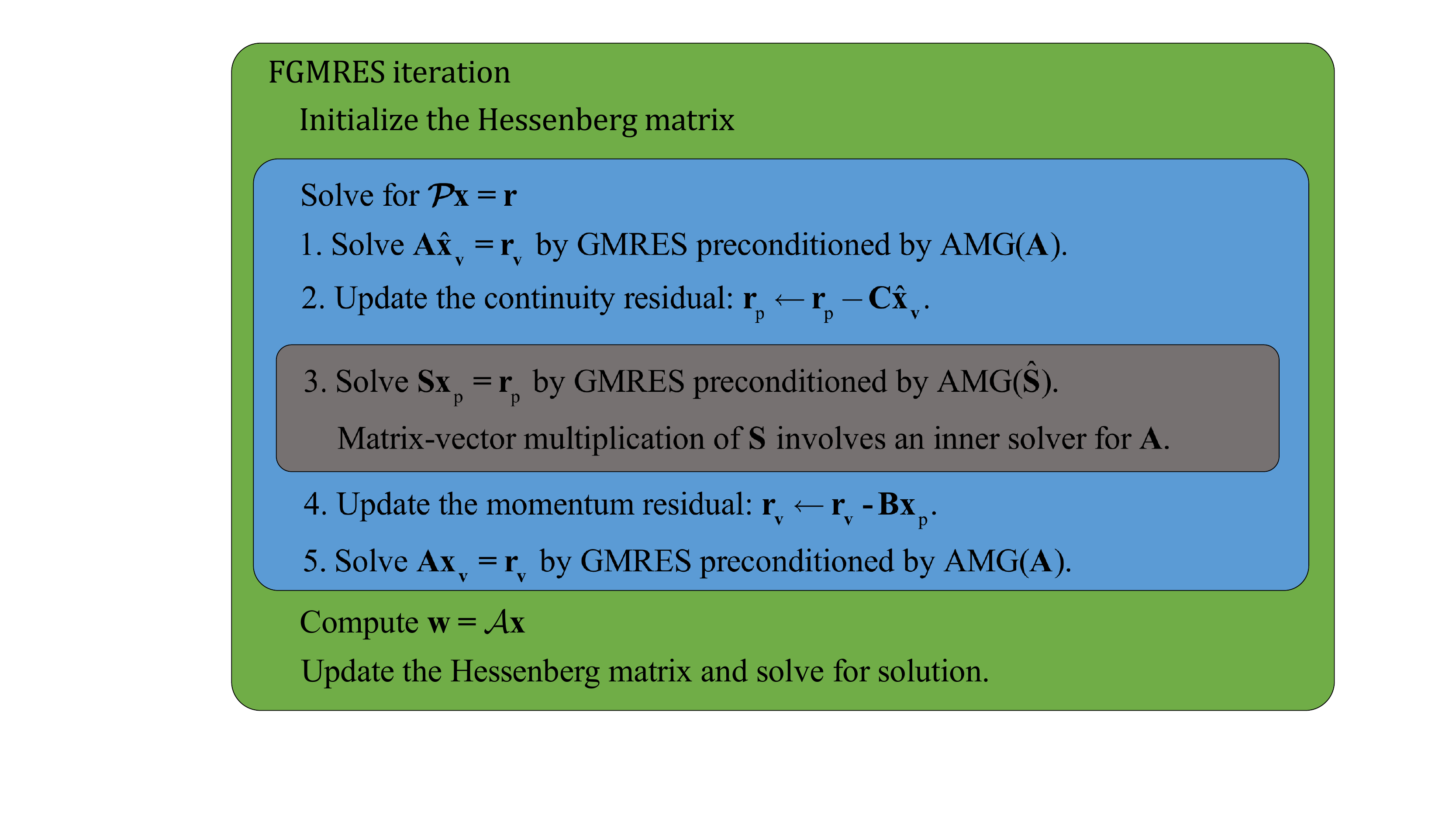} 
\caption{Implementation of the FGMRES with the nested block preconditioner. The green color represents the outer solver; the blue color represents the intermediate solver; the grey color represents the inner solver.} 
\label{fig:solver_flowchart}
\end{center}
\end{figure}

The FGMRES iteration for $\mathcal A$ serves as the \textit{outer solver} which tries to minimize the residual of $\mathcal A \bm x = \bm r$. Inside this FGMRES iteration, the application of $\hat{\mathcal P}_{SCR}$ is achieved through Algorithm \ref{algorithm:exact_block_factorization}, and one needs to solve with the block matrices $\boldsymbol{\mathrm A}$ and $\boldsymbol{\mathrm S}$ at this stage. We call it the \textit{intermediate solver}. When solving with the Schur complement, its action on a vector is defined by Algorihtm \ref{algorithm:matrix_free_mat_vec_for_S}, which necessitates using the \textit{inner solver} to solve with $\boldsymbol{\mathrm A}$. The three levels of solvers are illustrated in Figure \ref{fig:solver_flowchart} with different colors. 

\begin{remark}
In the construction of the proposed block preconditioners, the full $\mathcal L \mathcal D \mathcal U$ factorization of $\mathcal A$ is utilized. One can surely use only part of the factorization to devise different preconditioners. For example, the diagonal part $\mathcal D$ is an efficient candidate for the Stokes equations \cite{Silvester1994,Wathen1993}. Assuming exact arithmetic, it gives convergence within 4 iterations. Using the upper triangular part $\mathcal D \mathcal U$ often gives a good balance between the convergence rate and the computational cost \cite{Elman2002}, as it leads to convergence within 2 iterations \cite{Ipsen2001,Murphy2000}, assuming exact arithmetic. In our case, the full $\mathcal L \mathcal D \mathcal U$ block factorization gives the fastest convergence rate. We prefer this because the solution of the Schur complement equation is often the most expensive part of the overall algorithm. Therefore, in comparison with an upper triangular block preconditioner, we pay the price of solving the matrix problem $\boldsymbol{\mathrm A}$ twice with the purpose of mitigating the number of the solution procedure for the Schur complement.
\end{remark}

\begin{remark}
In the above algorithm, the nested block preconditioner $\hat{\mathcal P}_{SCR}$ can be regarded as a result of an inexact factorization of $\mathcal A$. The inexactness is due to the approximation made by the solvers in the intermediate and inner levels. The preconditioner is thus defined by the tolerances of these solvers. Using strict tolerances apparently makes $\hat{\mathcal P}_{SCR}$ closer to $\mathcal A$. However, this is impractical since this makes the algorithm as expensive as the SCR approach. On the other extreme, one may solve \eqref{eq:seg_sol_pres} by applying the preconditioner $\hat{\boldsymbol{\mathrm S}}$ once without invoking the inner solver. This makes the algorithm as simple as the coupled approach with the SIMPLE preconditioner and potentially endangers the robustness. We adjust the tolerances to tune the preconditioner, noting there is a lot of leeway in the choice of the tolerance value ranging from strict to loose. The effect of the tolerances of the intermediate and inner solvers will be studied in Section \ref{sec:numerical_results}.
\end{remark}

\begin{remark}
Choosing a good preconditioner for the Schur complement is critical for the performance of the proposed nested block preconditioner. In our experience, using a scaled pressure mass matrix gives satisfactory results as well \cite{Maliki2010}. For compressible materials, this preconditioner does not need to be explicitly assembled, and one can use $\boldsymbol{\mathrm D}$ directly (See \ref{app:consistent_linearization}). In this work, we focus on $\boldsymbol{\mathrm D} - \boldsymbol{\mathrm C} \left(\textup{diag}\left(\boldsymbol{\mathrm A}\right)\right)^{-1} \boldsymbol{\mathrm B}$, since this choice apparently is a better approximation of $\boldsymbol{\mathrm S}$. In \cite{Gurev2015}, a sparse approximate inverse is utilized to construct the preconditioner for the Schur complement, which is worth of future study.
\end{remark}

\section{Numerical Results}
\label{sec:numerical_results}
In our work, the outer solver is FGMRES($200$) with $n^{max} = 200$ and $\delta^a = 10^{-50}$. In the intermediate level, \eqref{eq:seg_sol_int_disp} and \eqref{eq:seg_sol_disp} are solved by GMRES($500$) preconditioned by AMG with $n^{max}_A = 500$ and $\delta^a_A = 10^{-50}$. The equation \eqref{eq:seg_sol_pres} is solved by GMRES($200$), with $n^{max}_S = 200$ and $\delta^a_S = 10^{-50}$. We use the AMG preconditioner constructed from $\hat{\boldsymbol{\mathrm S}}$. In the inner level, the linear system is solved via GMRES($500$) preconditioned by AMG with $n^{max}_I =500$ and $\delta^a_I = 10^{-50}$.  We use the BoomerAMG \cite{Henson2002} from the Hypre package \cite{Falgout2002} as the parallel AMG implementation. The settings of the BoomerAMG are summarized in Table \ref{table:amg_settings}. With the above settings, the accuracy of the solution is dictated by $\delta^r$, and the convergence rate is controlled by the tolerances $\delta^r_A$, $\delta^r_S$, and $\delta^r_I$.

\begin{table}[htbp]
\begin{center}
\renewcommand{\arraystretch}{1.0}
\begin{tabular}{p{9.5cm} p{5.5cm} }
\hline
Cycle type &  V-cycle \\
Coarsening method & HMIS \\
Interpolation method & Extended method (ext+i) \\
Truncation factor for the interpolation & $0.3$ \\
Threshold for being strongly connected & $0.5$ \\
Maximum number of elements per row for interp. & $5$ \\
The number of levels for aggressive coarsening & $2$ \\
\hline
\end{tabular}
\end{center}
\caption{Settings of the BoomerAMG preconditioner.}
\label{table:amg_settings}
\end{table}

To provide baseline examples, we solve the system of equations \eqref{eq:sym_diag_scaled_linear_system} by two different preconditioners. As the first example, we solve the the system of equations by FGMRES($200$) using $\hat{\mathcal P}_{\textup{SIMPLE}}$ with $n^{max} = 200$ and $\delta^a = 10^{-50}$. In this preconditioner, the settings of the linear solver (including the Krylov subspace method, the preconditioners, and the stopping criteria) associated with $\boldsymbol{\mathrm A}$ and $\hat{\boldsymbol{\mathrm S}}$ are exactly the same as the ones used in the nested block preconditioner. The accuracy of the solver is determined by $\delta^r$, and the performance of the preconditioner is controlled by $\delta^r_A$ and $\delta^r_S$. Notice that, in this preconditioner, $\delta^r_S$ is the tolerance for solving with the matrix $\hat{\boldsymbol{\mathrm S}}$.

As another baseline example, we choose to solve the linear system by GMRES($200$) preconditioned by a one-level additive Schwarz domain decomposition preconditioner \cite{Smith2004}. The maximum number of iterations is fixed at $10000$, and the tolerance for the absolute error is fixed at $10^{-50}$. In this preconditioner, each processor is assigned with a single subdomain, and an incomplete LU factorization (ILU) with a fill-in ratio 1.0 is invoked to solve the problem on the subdomains. This preconditioner is purely algebraic and is usually very competitive for medium-size parallel simulations. However, as will be shown in the numerical examples, the one-level domain decomposition preconditioner is not a robust option. Also, as the problem size and the number of subdomains grows, more iterations are needed to propagate information across the whole domain. In our implementation, the restricted additive Schwarz method from PETSc \cite{petsc-user-ref} is utilized as the domain decomposition preconditioner; the PILUT routine from Hypre \cite{Falgout2002} is used as the solver for the subdomain algebraic problem.

All numerical simulations are performed on the Stampede2 supercomputer at Texas Advanced Computing Center (TACC), using the Intel Xeon Platinum 8160 node. Each node contains 48 cores, with 2.1GHz nominal clock rate and 192GB RAM per node (4 GB RAM per core). 

\begin{figure}
	\begin{center}
	\begin{tabular}{cc}
\includegraphics[angle=0, trim=50 20 50 30, clip=true, scale = 0.35]{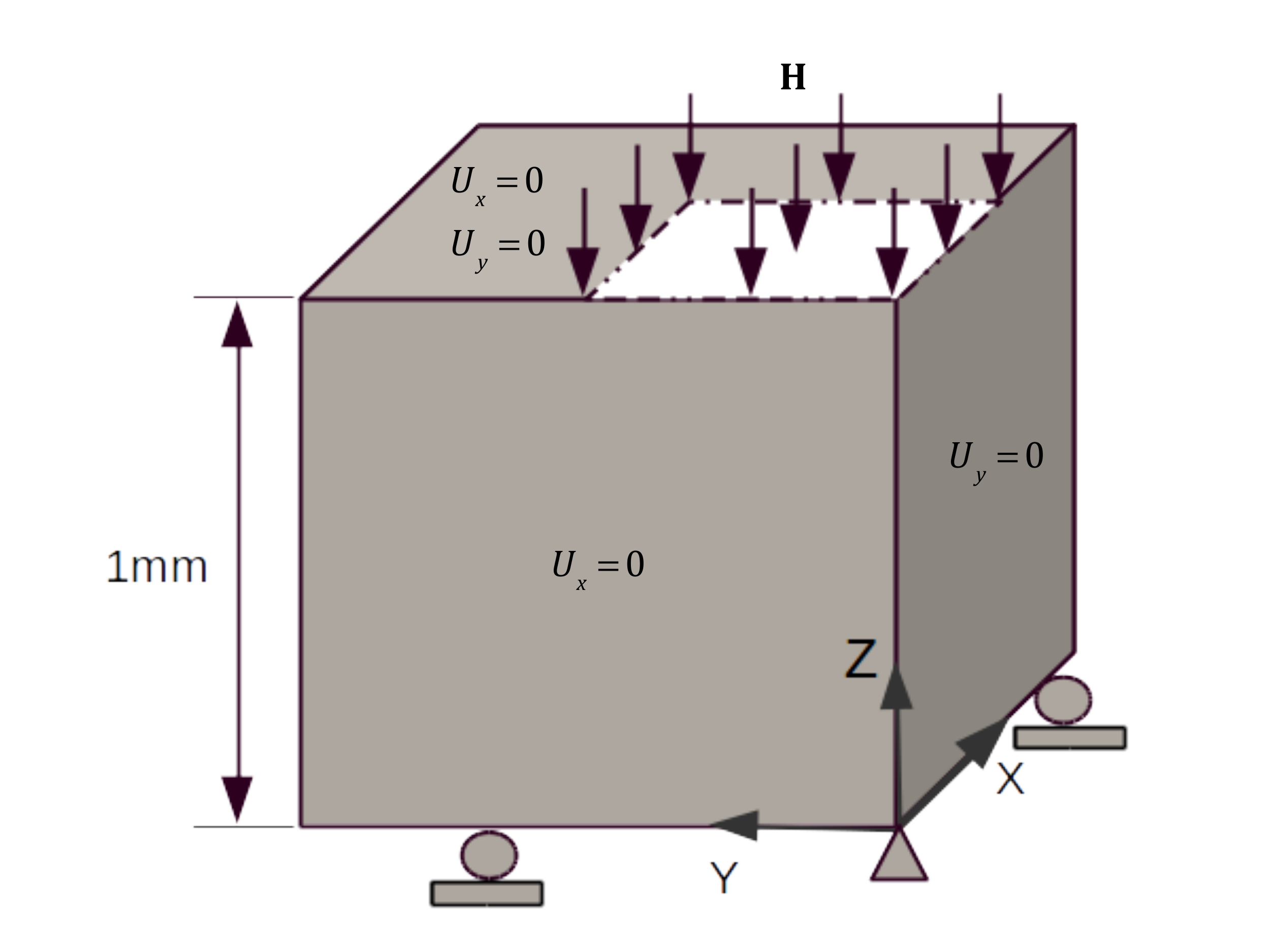} & 
\includegraphics[angle=0, trim=350 50 100 50, clip=true, scale = 0.18]{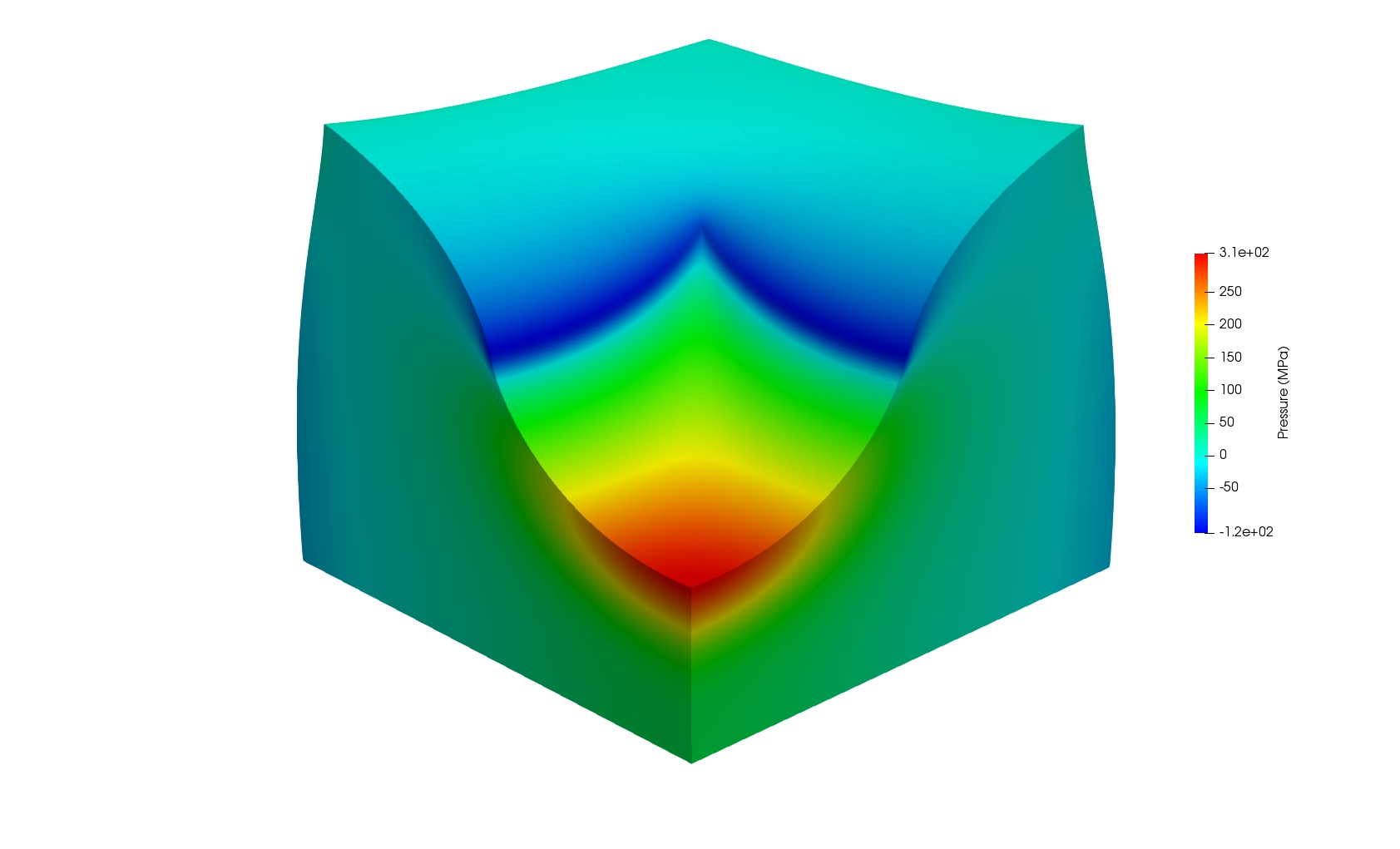}
\end{tabular}
\caption{Three-dimensional compression of a block: (left) geometry of the referential configuration and the boundary conditions; (right) pressure profile in the current configuration with $\Delta x = 1 / 3840$. } 
\label{fig:block_compression_setting}
\end{center}
\end{figure}

\subsection{Compression of a block}
\label{subsec:block_compression}
The compression of a unit block was proposed as a benchmark problem for nearly incompressible solids \cite{Reese2000}. The geometrical configuration and the boundary conditions are illustrated in Figure \ref{fig:block_compression_setting}. The problem is discretized in space by a uniform structured tetrahedral mesh generated by Gmsh \cite{Geuzaine2009}, and we use $\Delta x$ to denote the edge length of the mesh. The original benchmark problem was proposed in the quasi-static setting, and a `dead' surface load $\bm H$  is applied on a quarter portion of the top surface, pointing in the negative $z$-direction with magnitude $|\bm H|=320$ MPa. In this work, the problem is investigated in the dynamic setting by gradually increasing the load force as a linear function of time. The material is described by a Neo-Hookean model, whose Gibbs free energy function takes the form
\begin{align*}
G\left( \tilde{\bm C}, p \right) = \frac{\mu}{2\rho_0}  \left( \textup{tr}\tilde{\bm C} - 3 \right)  + \frac{p \sqrt{p^2+\kappa^2}-p^2}{2\kappa \rho_0} - \frac{\kappa}{2\rho_0}\ln \left( \frac{\sqrt{p^2+\kappa^2}-p}{\kappa} \right).
\end{align*}
Following \cite{Reese2000}, the material parameters are chosen as $\mu = 80.194$ MPa, $\kappa = 400889.806$ MPa, and $\rho_0 = 1.0 \times 10^3$ $\textup{kg/m}^3$. The corresponding Poisson's ratio $\nu$ is 0.4999. In Section \ref{subsec:solver_different_material_para}, we examine the robustness of the preconditioner with regard to varying material moduli. In the following discussion, the governing equations have been non-dimensionalized by the centimetre-gram-second units. Note that the edge length of the cube is $1$ mm $=0.1$ cm. Then the number of elements in each direction of the cube is given by $1/(10\Delta x)$.

\subsubsection{Performance with varying inner solver accuracy}
\label{subsec:inner_solver}
In this test, we investigate the impact of the accuracy of the inner solver on the overall iterative method. We fix the mesh size to be $\Delta x = 1/640$ and the time step size to be $\Delta t = 10^{-1}$. The simulation is performed with 8 CPUs, with approximately 131072 equations assigned to each CPU.  In this study, we choose $\delta^r = 10^{-8}$, and we consider two settings for the intermediate solver: $\delta^r_A = \delta^r_S = 10^{-10}$ and $\delta^r_A = \delta^r_S = 10^{-6}$. We collect the statistics of the solver in the first time step with varying values of $\delta^r_I$ (Table \ref{table:impact_inner_on_FGMRES}). 
The results associated with $\delta^r_I = 10^0$ are obtained by solving $\hat{\boldsymbol{\mathrm S}} \bm x_p = \hat{\bm r}_{p}$ in step 3 of Algorithm \ref{algorithm:exact_block_factorization}. 

In our numerical experiments, we observe that with the choice of $\delta^r_A = \delta^r_S = \delta^r_I = 10^{-10}$, the outer solver converges in less than two iterations on average. In fact, we also experimented with stricter tolerances and observed convergence of the outer solver in one iteration. (We do not report this because this stricter choice requires larger size of the Krylov subspace which is incompatible with our current settings.) This result corroborates the fact that the full $\mathcal L \mathcal D \mathcal U$ block preconditioner gives convergence in one iteration with exact arithmetic.

In the literature, the choice for the inner solver accuracy is under debate. In \cite{May2008}, it is suggested that the inner solver should be more accurate than its upper-level counterpart (i.e., $\delta^r_I \leq \delta^r_S$ in our case) to guarantee accurate representation of the Schur complement. Meanwhile, it is shown in \cite{Bouras2005} that the Krylov methods are in fact very robust under the presence of inexact matrix-vector multiplications. In our test, as we gradually release the tolerance $\delta^r_I$, it is observed that the inner solver converges with fewer iterations while the outer solver requires more iterations to reach convergence to compensate for the inaccurate evaluations of the Schur complement. As $\delta^r_I$ gets larger than $\delta^r_S$, initially the overhead is low. As the tolerance further increases, the outer solver requires more iterations and the overall cost of the solver grows correspondingly. For the two cases, the break-even points are achieved with $\delta^r_I = 10^{-6}$ and $10^{-4}$, respectively. Examining the number of iterations for the outer solver, we observe a steady growth of $n$ once $\delta^r_I$ grows larger than $\delta^r_S$. Though it is hard to predict the optimal choice of $\delta^r_I$ for general cases, we observe that a choice of $\delta^r_I = \delta^r_S$ is safe for robust performances; a slightly relaxed tolerance for the inner solver (e.g. $\delta^r_I = 10^2 \delta^r_S$) is beneficial for efficiency.

For comparison, we also examined the solver performance without the inner solver. We solve with $\hat{\boldsymbol{\mathrm S}}$ instead of $\boldsymbol{\mathrm S}$ in \eqref{eq:seg_sol_pres} directly. This corresponds to a highly inaccurate evaluation of the Schur complement. We see that the iteration number and the CPU time of the outer solver both increase significantly. The severe degradation of solver performance signifies the importance of an accurate evaluation of the Schur complement.

\begin{table}[htbp]
\begin{center}
\tabcolsep=0.19cm
\renewcommand{\arraystretch}{1.2}
\begin{tabular}{P{3.0cm} P{1.0cm} P{3.0cm} P{0.8cm} P{0.8cm} P{1.0cm} P{1.0cm} P{1.0cm} }
\hline
& $\delta^r_I$ & CPU time (sec.) & $\hat{l}$ &  $n$ & $\bar{n}_A$ & $\bar{n}_S$ & $\bar{n}_I$   \\
\hline
$\delta^r_A = \delta^r_S = 10^{-10}$ & $10^{0}$ & $4.86 \times 10^3$ & 4 & 477 & 74.52 & 33.89 & -  \\
& $10^{-2}$ & $9.02 \times 10^2$ & 4 & 17 & 75.62 & 22.29 & 29.31  \\
& $10^{-4}$ & $8.08 \times 10^2$ & 4 & 11 & 75.30 & 22.27 & 45.00 \\
& $10^{-6}$ & $6.97 \times 10^2$ & 4 & 8 & 75.19 & 23.13 & 55.82  \\
& $10^{-8}$ & $8.11 \times 10^2$  & 4 & 8 & 75.19  & 22.13  & 65.15  \\
& $10^{-10}$ & $8.47 \times 10^2$ & 4 & 7 & 74.86 & 23.29 & 74.62  \\[0.5em]
$\delta^r_A = \delta^r_S = 10^{-6}$ & $10^{0}$ &  $4.87 \times 10^3$  & 4 & 664 & 55.60 & 21.29 & -  \\
& $10^{-2}$ & $6.30 \times 10^2$  & 4 & 18 & 56.68 & 13.74 & 30.35  \\
& $10^{-4}$ & $5.10 \times 10^2$ & 4 & 11 & 56.30 & 13.73 & 46.14 \\
& $10^{-6}$ & $5.12 \times 10^2$  & 4 & 9 & 56.06 & 14.11 & 56.91 \\
& $10^{-8}$ & $6.01 \times 10^2$ & 4  & 9 & 56.17 & 14.44  & 65.62  \\
& $10^{-10}$ & $6.82 \times 10^2$ & 4 & 9 & 56.17 & 14.56 & 74.87  \\
\hline
\end{tabular}
\end{center}
\caption{The impact of the accuracy of the inner solver on the performance of the linear solver. The CPU time is collected for the linear solver only; $\hat{l}$ represents the total number of nonlinear iterations; $n$ represents the total number of FGMRES iterations; $\bar{n}_A$ represents the averaged number of iterations for solving with $\boldsymbol{\mathrm A}$ in \eqref{eq:seg_sol_int_disp} and \eqref{eq:seg_sol_disp}; $\bar{n}_S$ represents the averaged number of iterations for solving \eqref{eq:seg_sol_pres}; $\bar{n}_I$ represents the averaged number of iterations for solving \eqref{eq:S_inner_A_eqn}.}
\label{table:impact_inner_on_FGMRES}
\end{table}

\subsubsection{Performance with varying intermediate solver accuracy}
\label{subsec:conv_rate_isotropic_compression}
\begin{figure}
	\begin{center}
	\begin{tabular}{cc}
\includegraphics[angle=0, trim=80 100 360 280, clip=true, scale = 0.04]{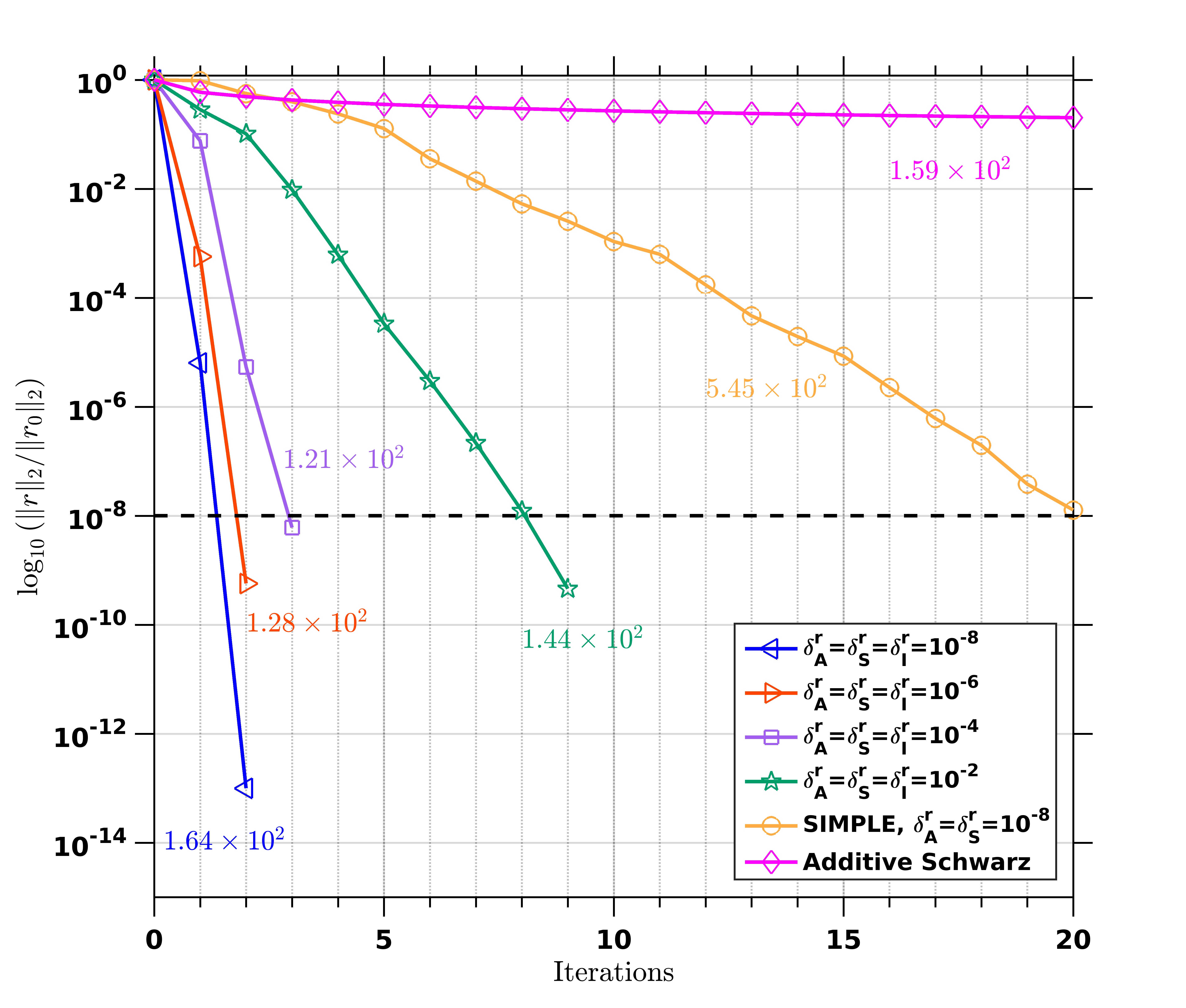} &
\includegraphics[angle=0, trim=80 100 360 280, clip=true, scale = 0.04]{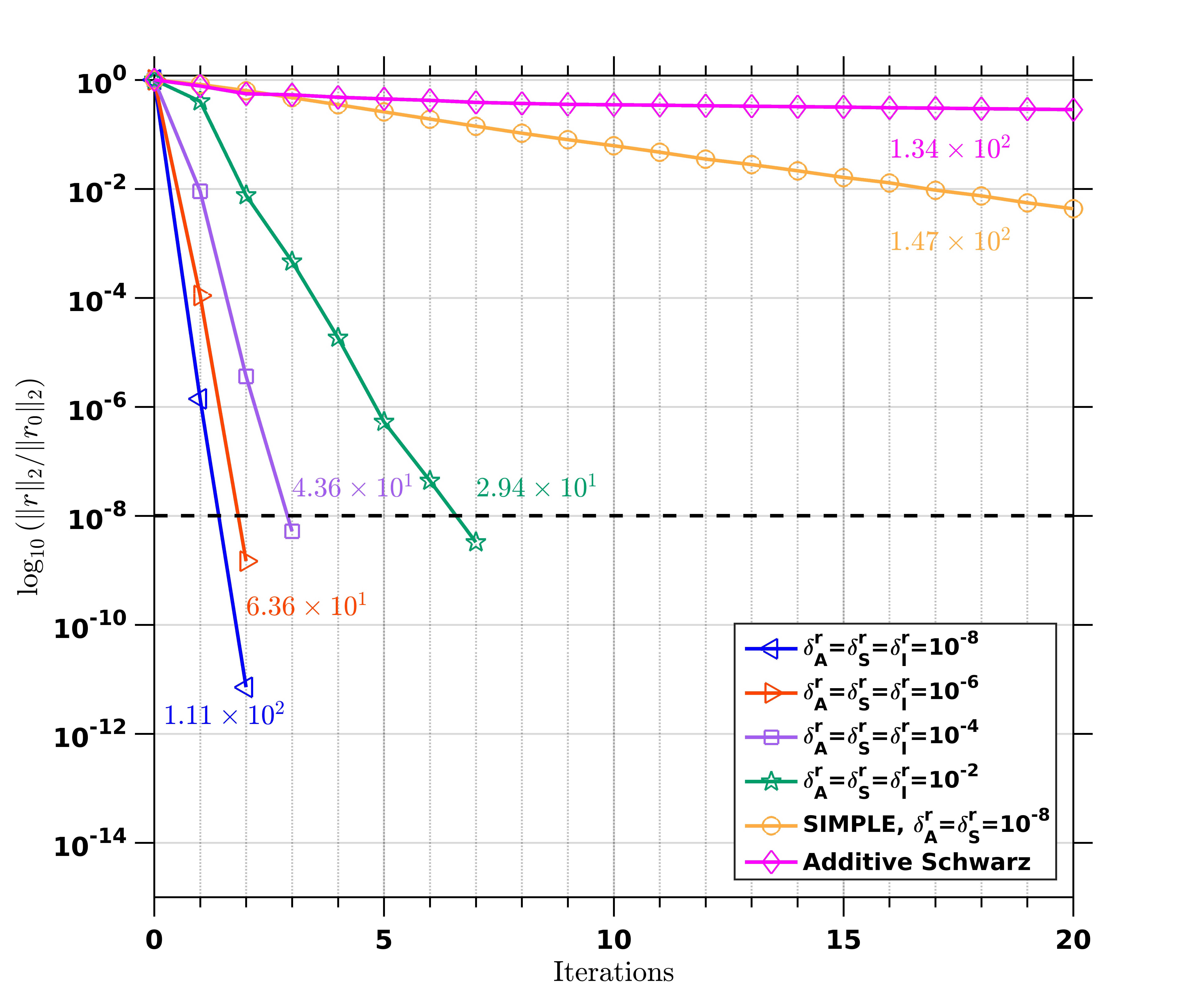}
\end{tabular}
\caption{Convergence history for $\Delta t = 10^{-1}$ (left) and $10^{-5}$ (right). The horizontal dashed black line indicates the prescribed stopping criterion for the relative error, which is $10^{-8}$ here. In the case of $\Delta t = 10^{-1}$, the SIMPLE method converges in 21 iterations, and the additive Schwarz method converges in 2644 iterations. In the case of $\Delta t = 10^{-5}$, the SIMPLE method converges in 71 iterations, and the additive Schwarz method converges in 2030 iterations. The numbers indicate the averaged time per nonlinear iteration in seconds.} 
\label{fig:conv_history_m64}
\end{center}
\end{figure}

In this example, we examine the effect of varying intermediate solver tolerances on the solver performance. We consider a uniform mesh with $\Delta x = 1/640$, with two time step sizes: $\Delta t = 10^{-1}$ and $10^{-5}$. We choose $\delta^r = 10^{-8}$ for the outer solver. We set $\delta^r_A = \delta^r_S = \delta^r_I$ and vary their values from $10^{-8}$ to $10^{-2}$. To make comparisons, the same problem is simulated with the SIMPLE preconditioner and the additive Schwarz preconditioner. In the SIMPLE preconditioner, we solve the equations associated with $\boldsymbol{\mathrm A}$ and $\hat{\boldsymbol{\mathrm S}}$ with $\delta^r_A = \delta^r_S = 10^{-8}$. The convergence is monitored for the first time step, which is usually the most challenging part of dynamic calculations. The convergence history of the linear solver in the first nonlinear iteration is plotted in Figure \ref{fig:conv_history_m64}. It can be seen that the accuracy of the intermediate solvers affects the convergence rate of the linear solver. When the equations in the intermediate level are solved to a high precision, the convergence rate of the outer solver is steep. As one looses the tolerances for the intermediate solvers, the proposed algorithm requires more iterations for convergence. Yet, even for the tolerance as loose as $10^{-2}$, the convergence rate is still much steeper than that of the SIMPLE preconditioner. The average time for solving the matrix problem per nonlinear iteration is reported in the figures as well. We observe that when choosing a strict tolerance for the intermediate and inner solvers, although convergence is achieved with fewer iterations, the cost per iteration is high and the overall time to solution is correspondingly high. A looser tolerance renders the application of the nested block preconditioner more cost-effective, and the overall algorithm is faster. In comparison with the SIMPLE and additive Schwarz methods, the proposed nested block preconditioning technique is fairly competitive.

\subsubsection{Performance with varying material properties}
\label{subsec:solver_different_material_para}
In this example, we vary the material properties and study the robustness of the proposed preconditioner. The Poisson's ratio $\nu$ varies from $0.0$ to $0.5$, spanning the range relevant to most engineering and biological materials. The shear modulus $\mu$ is taken as $80.194 \times \eta$ MPa, wherein $\eta$ is a non-dimensional number. Correspondingly, the compression force is adjusted by multiplying with the scaling factor $\eta$ for values of $10^{-2}$, $10^0$, and $10^2$. The stopping condition for the linear solver is $\delta^r = 10^{-8}$, and we choose $\delta^r_A = \delta^r_S = \delta^r_I = 10^{-6}$. The mesh size is fixed to be $\Delta x = 1/480$, and the problem is simulated with $8$ CPUs. The time step size is $\Delta t = 10^{-1}$ and we integrate the problem up to $T=1.0$. We use a relatively large time step size here to make the matrix $\boldsymbol{\mathrm A}$ dominated by the stiffness matrix. The statistics of the solver performance are collected over ten time steps. The averaged number of iterations as well as the averaged CPU time for one nonlinear iteration $\bar{T}_L$ are reported in Table \ref{table:material_nu_E}.

\begin{table}[htbp]
\begin{center}
\tabcolsep=0.19cm
\renewcommand{\arraystretch}{1.2}
\begin{tabular}{P{2.8cm} P{4cm} P{4cm} P{4cm}}
\hline
$\bar{n}$ [$\bar{n}_A$, $\bar{n}_S$] ($\bar{T}_L$) & $\eta=10^{-2}$ & $\eta=10^{0}$ & $\eta=10^2$ \\
 $\nu = 0.0$ & 2.0 [46.9, 15.9] (46.6) & 2.0 [48.1, 16.0] (48.3) & 2.0 [47.9, 15.3] (46.2) \\
 $\nu = 0.1$ & 2.0 [48.5, 19.0] (49.2) & 2.0 [48.4, 17.9] (50.3) & 2.0 [48.1, 15.5] (46.5) \\
 $\nu = 0.2$ & 2.0 [48.3, 20.2] (52.8) & 2.0 [48.0, 19.9] (52.9) & 2.0 [48.3, 16.8] (47.0) \\
 $\nu = 0.3$ & 2.0 [47.9, 23.1] (56.4) & 2.0 [41.1, 21.5] (57.9) & 2.0 [48.5, 17.7] (48.6) \\
 $\nu = 0.4$ & 2.0 [47.4, 28.5] (66.4) & 2.0 [48.2, 25.8] (65.5) & 2.0 [48.6, 19.2] (50.6) \\
 $\nu = 0.5$ & 2.2 [47.1, 36.3] (101.2) & 2.0 [47.4, 24.6] (66.3) & 2.0 [46.5, 20.3] (48.6)  \\
\hline
\end{tabular}
\end{center}
\caption{The performance of the linear solver with varying material properties.}
\label{table:material_nu_E}
\end{table}

For all cases, the number of iterations for the outer solver maintains around two. In fact, it is only for the case of $\nu=0.5$ and $\eta = 10^{-2}$ that the outer solver needs slightly more than two iterations. The number of iterations for solving with $\boldsymbol{\mathrm A}$ in \eqref{eq:seg_sol_int_disp} and \eqref{eq:seg_sol_disp} is maintained around $47$, and hence can be regarded as independent with respect to the material property. The number of iterations for solving \eqref{eq:seg_sol_pres} increases with increasing the Poisson's ratio. This can be explained by looking at $\boldsymbol{\mathrm S} = \boldsymbol{\mathrm D} - \boldsymbol{\mathrm C} \boldsymbol{\mathrm A}^{-1} \boldsymbol{\mathrm B}$. The matrix $\boldsymbol{\mathrm D}$ is dominated by the mass matrix scaled with a factor of $\beta$. As $\nu$ approaches $0.5$, the isothermal compressibility coefficient $\beta$ goes to zero. Consequently, the well-conditioned matrix $\boldsymbol{\mathrm D}$ diminishes, and the condition number of the Schur complement gets larger. This is reflected in the increase of $\bar{n}_S$ as $\nu$ goes from $0.0$ to $0.5$ for all three shear moduli. On the other hand, $\bar{n}_S$ increases as the material gets softer, and this trend is pronounced as the Poisson's ratio gets larger. This can be explained by looking at $\boldsymbol{\mathrm A}^{-1}$ in the Schur complement. For large time steps, $\boldsymbol{\mathrm A}$ contains a significant contribution from the stiffness matrix, and the inverse of the stiffness matrix is proportional to $1/\mu$. It is known that $\textup{diag}\left(\boldsymbol{\mathrm A}\right)$ is not a good candidate for approximating the stiffness matrix, and this is magnified for softer materials due to the factor $1/\mu$.

\subsubsection{Parallel performance}
\label{subsec:parallel_performance}
We investigate the efficiency of the method by evaluating the fixed-size scalability performance. The spatial mesh size is $\Delta x = 1/1280$, with about $8.39 \times 10^{6}$ degrees of freedom. The time step size is fixed at $10^{-5}$, and we integrate the problem in time up to $T=10^{-4}$. The stopping criterion for the FGMRES iteration is $\delta^r = 10^{-3}$; the tolerances for the intermediate and inner solvers are $\delta^r_A = \delta^r_S =\delta^r_I = 10^{-3}$. We observe that the efficiency of the numerical simulation is maintained at a high level (around $90\%$) for a wide range of processor counts (Table \ref{table:strong_scaling}).
\begin{table}[htbp]
\begin{center}
\tabcolsep=0.19cm
\renewcommand{\arraystretch}{1.2}
\begin{tabular}{P{2.0cm} P{2.0cm} P{3.0cm} P{2.5cm} P{2.0cm} }
\hline
Proc. & $T_A$ (sec.) & $T_L$ (sec.) & Total (sec.) & Efficiency   \\
\hline
2 & $3.13\times 10^3$ & $2.16 \times 10^4$ & $2.49\times 10^4$ & $100\%$ \\
4 & $1.57\times 10^3$ & $1.09 \times 10^4$ & $1.26\times 10^4$ & $99\%$ \\
8 & $8.49\times 10^2$ & $5.58 \times 10^3$ & $6.48\times 10^3$ & $96\%$ \\
16 & $4.38\times 10^2$ & $2.96 \times 10^3$ & $3.43\times 10^3$ & $91\%$ \\
32 & $2.33\times 10^2$ & $1.62 \times 10^3$ & $1.87\times 10^3$ & $83\%$ \\
64 & $1.10\times 10^2$ & $8.37 \times 10^2$ & $9.56\times 10^2$ & $81\%$ \\
128 & $5.65\times 10^1$ & $3.84 \times 10^2$ & $4.49\times 10^2$ & $87\%$ \\
\hline
\end{tabular}
\end{center}
\caption{The strong scaling performance. $T_A$ and $T_L$ represent the timings for matrix assembly and linear solver, respectively. The efficiency is computed based on the total time.}
\label{table:strong_scaling}
\end{table}

\begin{table}
\begin{center}
\tabcolsep=0.19cm
\renewcommand{\arraystretch}{1.2}
\begin{tabular}{@{\extracolsep{4pt}}P{0.8cm} P{1.0cm} P{0.6cm} P{0.8cm} P{0.8cm} P{1.0cm} P{1.0cm} P{1.0cm} P{1.5cm} P{1.0cm}@{}}
\hline
\multirow{2}{*}{$\frac{1}{\Delta x}$} & \multirow{2}{*}{Proc.} & \multicolumn{4}{c}{$\hat{\mathcal P}_{SCR}$} & \multicolumn{2}{c}{SIMPLE} & \multicolumn{2}{c}{Additive Schwarz}\\
\cline{3-6} \cline{7-8}\cline{9-10}
& & $\bar{n}$ & $\bar{n}_A$ & $\bar{n}_S$ & $\bar{T}_L$ & $\bar{n}$ & $\bar{T}_L$ & $\bar{n}$ & $\bar{T}_L$
\\
\hline
\multicolumn{10}{l}{$\Delta t = 10^{-1}$} \\
$480$ & 8 & 2.3 & 31.7 & 6.7 & 19.7 & 13.3 & 21.6 & 1114.4 & 20.4 \\
$960$ & 64 & 2.5 & 43.1 & 7.4 & 50.2 & 17.9 & 63.6 & 3368.9 & 106.8 \\
$1920$ & 512 & 2.7 & 55.4 & 9.1 & 108.0 & 25.0 & 153.6 & 8642.4 & 305.4 \\
$3840$ & 4096 & 2.9 & 68.8 & 9.6 & 220.8 & 47.6 & 504.2 & NC & NC \\  
\hline
\multicolumn{10}{l}{$\Delta t = 10^{-5}$} \\
$480$ & 8 & 2.3 & 4.6 & 16.1 & 5.3 & 22.7 & 6.3 & 916.0 & 17.0 \\  
$960$ & 64 & 2.0 & 6.9 & 26.4 & 18.5 & 38.6 & 31.6 & 2133.9 & 67.4 \\ 
$1920$ & 512 & 2.0 & 9.1 & 34.3 & 52.8 & 65.7 & 71.0 & 9669.1 & 315.0 \\
$3840$ & 4096 & 2.2 & 11.3 & 42.0 & 139.0 & 101.2 & 221.4 & NC & NC \\
\hline
\end{tabular}
\end{center}
\caption{Comparison of the averaged iteration counts and CPU time in seconds for the nested block preconditioner $\hat{\mathcal P}_{SCR}$, the SIMPLE preconditioner, and the additive Schwarz preconditioner. NC stands for no convergence. For the $\Delta x = 1/3840$ case, the additive Schwarz preconditioner failed to converge in 10000 iterations.}
\label{table:weak_scaling}
\end{table}

To compare the performance of different preconditioners, we also perform a weak scaling test of the solver, with $\delta^r = 10^{-3}$. Tolerances are set to $\delta^r_A = \delta^r_S = \delta^r_I = 10^{-3}$ for the nested block preconditioner and to $\delta^r_A = \delta^r_S = 10^{-3}$ for the SIMPLE preconditioner. The computational mesh is progressively refined and each CPU is assigned approximately $5.53\times 10^4$ equations. We simulate the problem with two different time step sizes: $\Delta t = 10^{-1}$ and $10^{-5}$. The statistics of the solver performance are collected for ten time steps (Table \ref{table:weak_scaling}). We observe that the iteration counts for the outer solver using the nested block preconditioner are independent of mesh refinement. At large time steps, $\boldsymbol{\mathrm A}$ is dominated by the stiffness matrix and its solution procedure requires more iterations. In the meantime, the Schur complement has a better condition number and converges with fewer iterations. At small time steps, the situation is opposite. The matrix $\boldsymbol{\mathrm A}$ is dominated by the mass matrix, and it can be solved with fewer iterations. The mesh refinement has an impact on the intermediate solvers, and we observe an increase of the number of iterations in $\bar{n}_A$ and $\bar{n}_S$. For the SIMPLE preconditioner, the iteration counts and the CPU time grow faster than those of the nested block preconditioner. The additive Schwarz method converges faster per iteration. However, the number of iterations for convergence is much higher. For the finest mesh, the additive Schwarz method fails to converge in 10000 iterations. The proposed nested block preconditioner gives the most robust and efficient performance.

\subsection{Tensile test of an anisotropic fibre-reinforced hyperelastic soft tissue model}
\begin{figure}
	\begin{center}
	\begin{tabular}{cc}
\includegraphics[angle=0, trim=20 60 90 30, clip=true, scale = 0.35]{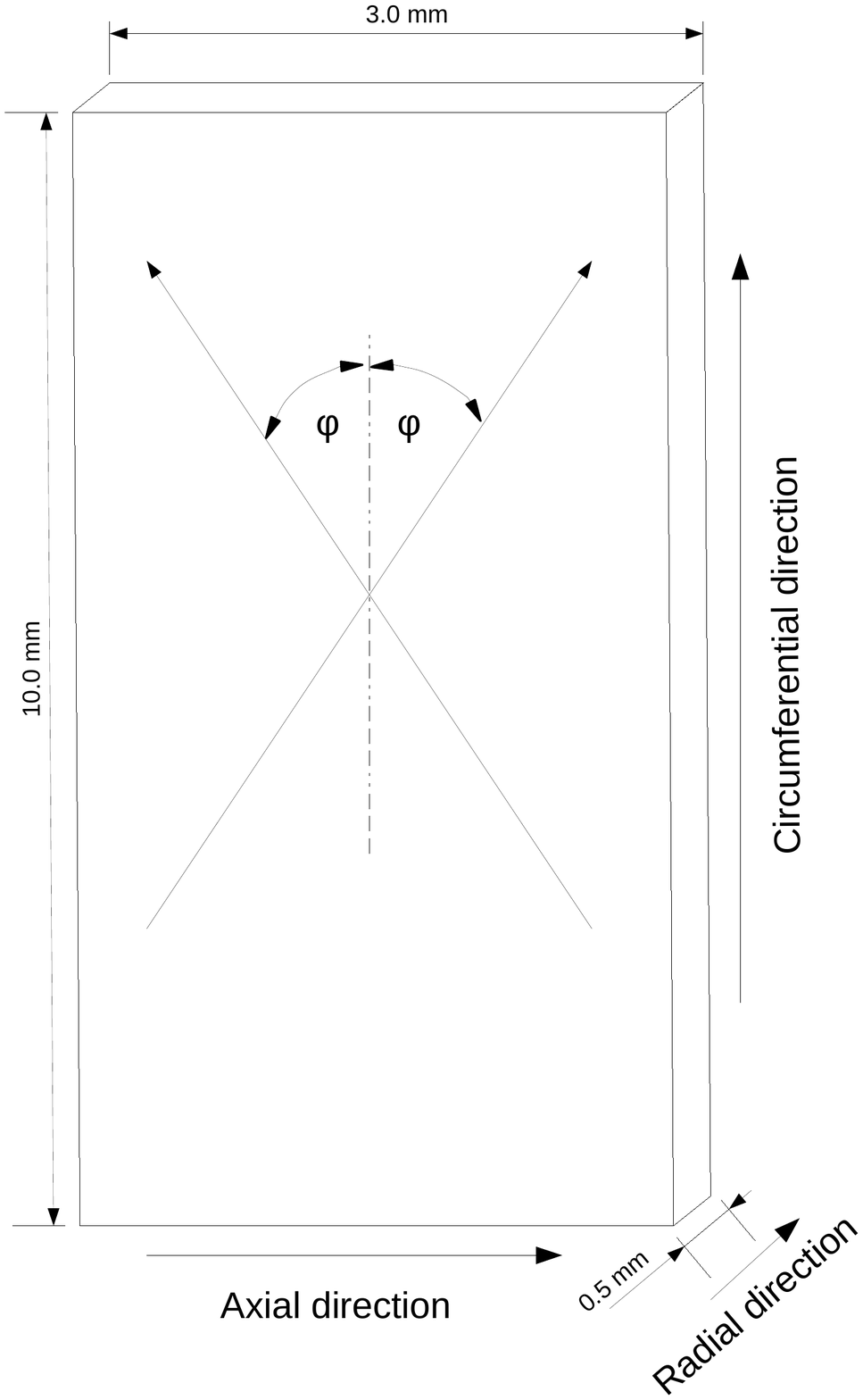} &
\includegraphics[angle=0, trim=100 0 100 0, clip=true, scale = 0.05]{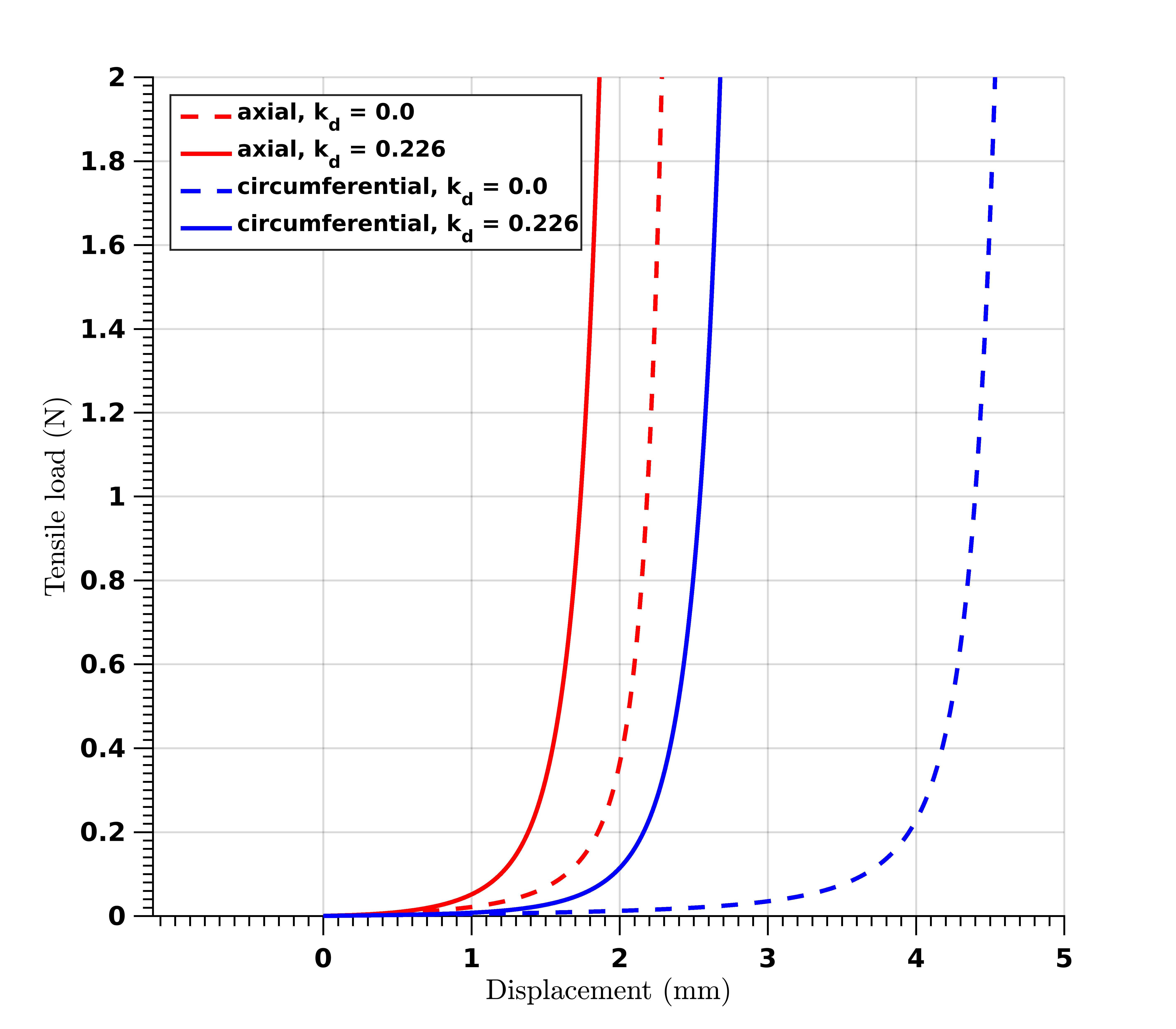} \\
(a) & (b)
\end{tabular}
\caption{Three-dimensional tensile test of an iliac adventitial strip: (a) geometry of the referential configuration; (b) computed load-displacement curves of the circumferential (blue) and axial specimens (red) with ($\kappa = 0.226$, solid curves) and without ($\kappa = 0.0$, dashed curves) dispersion of the collagen fibres.} 
\label{fig:tensile_test_setting}
\end{center}
\end{figure}

In this example, we apply the proposed preconditioning technique to an anisotropic hyperelastic material model, which has been used to describe arterial tissue layers with distributed collagen fibres. The isochoric and volumetric parts of the free energy are
\begin{align*}
G_{iso}(\tilde{\bm C} ) &= G_{iso}^{g}( \tilde{\bm C} ) + \sum_{i=1,2} G_{iso}^{f_i}( \tilde{\bm C} ), \quad G_{vol}(p) = \frac{p}{\rho_0}, \\
G_{iso}^{g}( \tilde{\bm C} ) &= \frac{\mu}{2\rho_0} \left( \textup{tr}\tilde{\bm C} -3 \right), \quad
G_{iso}^{f_i}( \tilde{\bm C} ) = \frac{k_1}{2k_2\rho_0} \left( e^{k_2 \bar{E}_i^2} - 1 \right), \\
\bar{E}_i &= \bm H_i : \tilde{\bm C} - 1, \quad \bm H_i = k_d \bm I + (1-3k_d)(\bm a_i \otimes \bm a_i).
\end{align*}
In the above, $G_{iso}^{g}$ models the groundmatrix via an isotropic Neo-Hookean material, with $\mu$ being the shear modulus; $G_{iso}^{f_i}$ models the $i$th family of collagen fibres by an exponential function. In $G_{iso}^{f_i}$, $\bm a_i$ is a unit vector that describes the mean orientation of the $i$th family of fibres in the reference configuration. The parameter $k_d \in [0, 1/3]$ is a structural parameter that characterizes the dispersion of the collagen fibres. For ideally aligned fibres, the dispersion parameter $k_d$ is 0, while for isotropically distributed fibres, it takes the value $1/3$. The parameter $k_1$ is a material parameter that describes the stiffness of the fibre, and $k_2$ is a non-dimensional parameter. The volumetric energy $G_{vol}$ indicates that the model is fully incompressible. Interested readers are referred to \cite{Gasser2006} for detailed discussions of the histology and constitutive modeling of the arterial layers. In the numerical study, we perform a tensile test for the tissue model. Following \cite{Gasser2006}, the geometry of the specimen has length $10.0$ mm, width $3.0$ mm, and thickness $0.5$ mm. The material parameters are $\mu=7.64$ kPa, $k_1=996.6$ kPa, $k_2 = 524.6$. Assuming that the fibre orientation has no radial component, the unit vector is characterized completely by $\varphi$, the angle between the circumferential direction and the mean fibre orientation direction (see Figure \ref{fig:tensile_test_setting} (a)). For the circumferential specimen, $\varphi = 49.98^\circ$; for the axial specimen, $\varphi = 40.02^\circ$. On the loading surface, traction force is applied and the face is constrained to move only in the loading direction. Symmetric boundary conditions are properly applied, and we only consider one-eighth of the specimen in the simulations. 

\begin{figure}
	\begin{center}
	\begin{tabular}{cccc}
\includegraphics[angle=0, trim=800 30 600 30, clip=true, scale = 0.23]{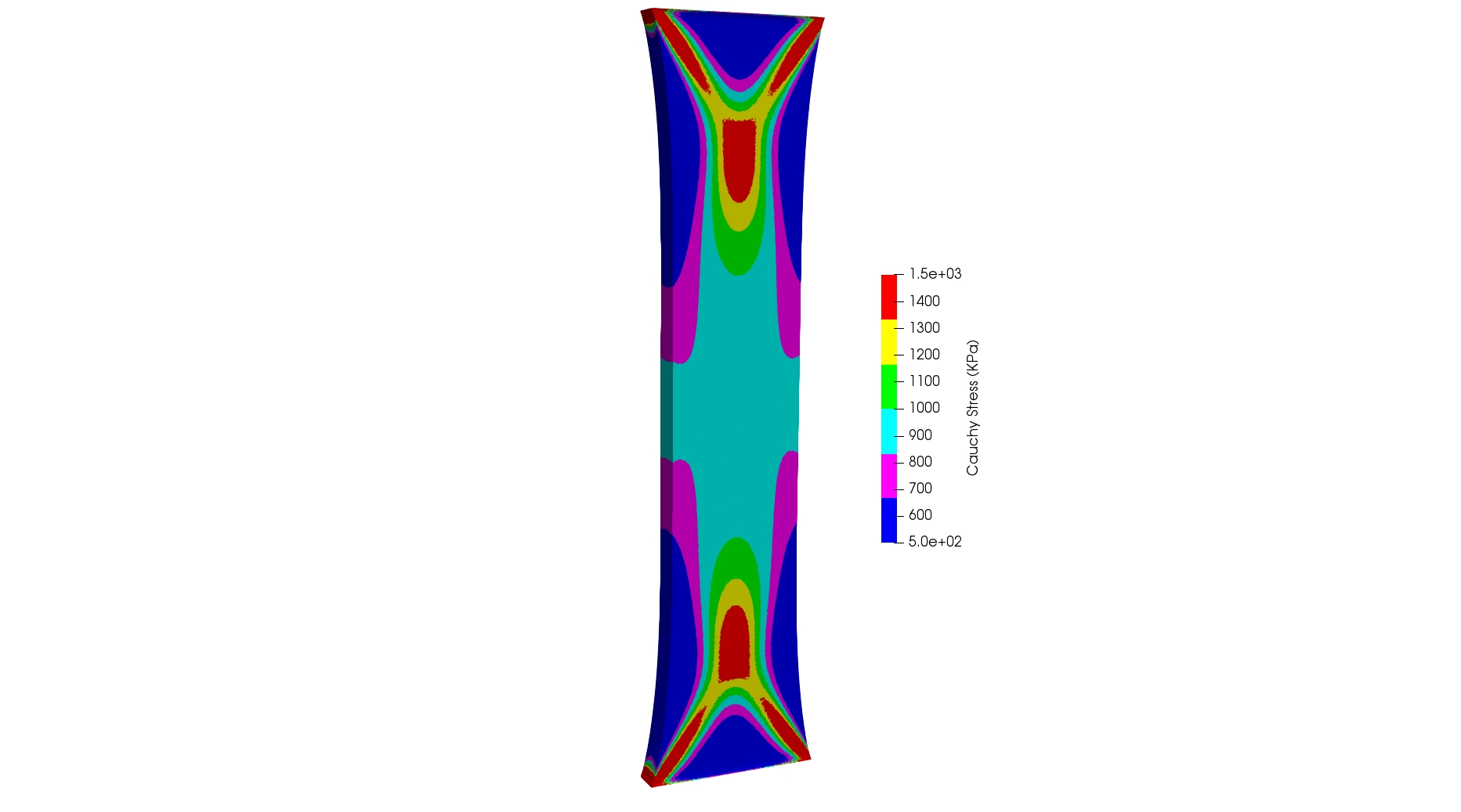} &
\includegraphics[angle=0, trim=800 30 700 30, clip=true, scale = 0.23]{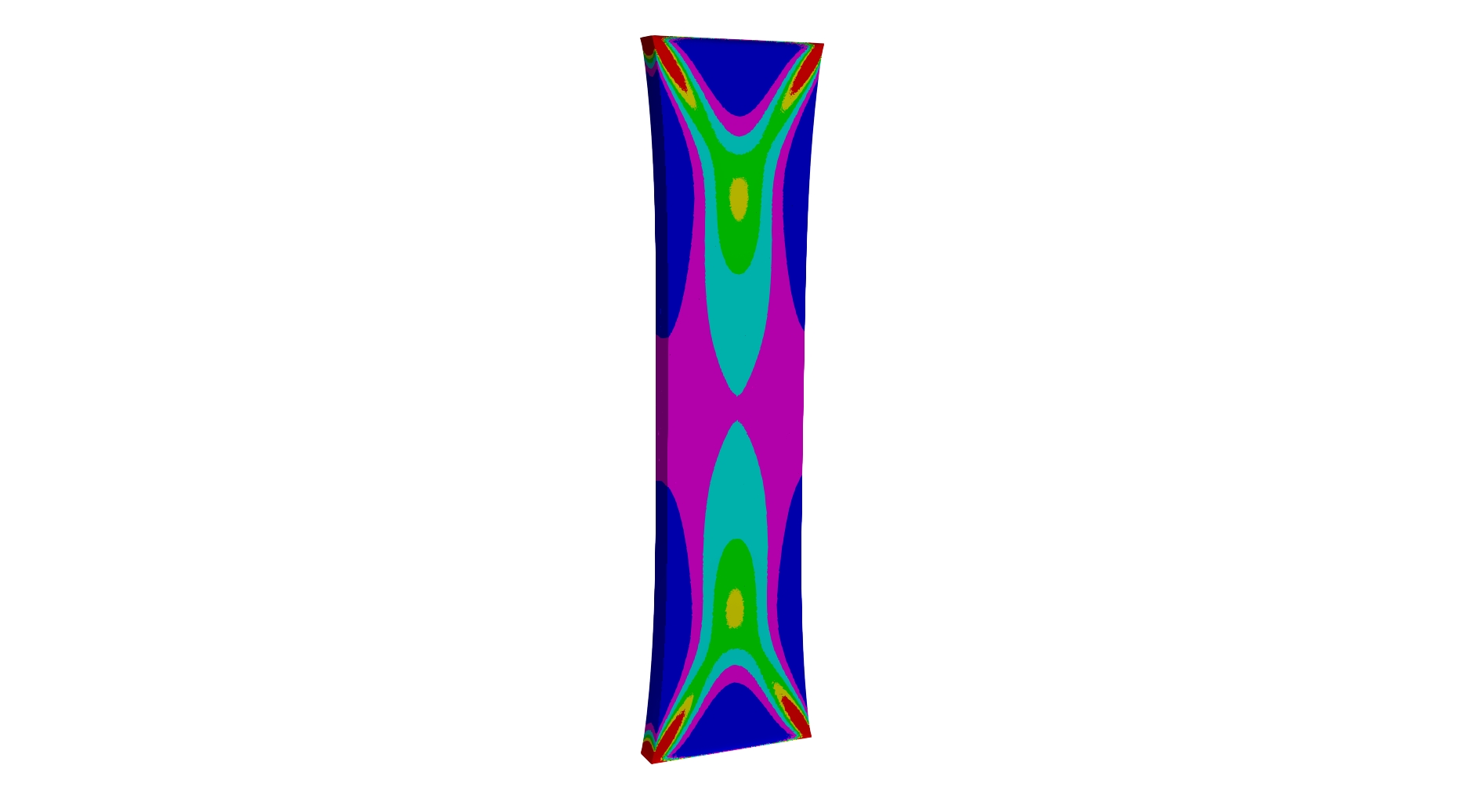} &
\includegraphics[angle=0, trim=800 30 600 30, clip=true, scale = 0.23]{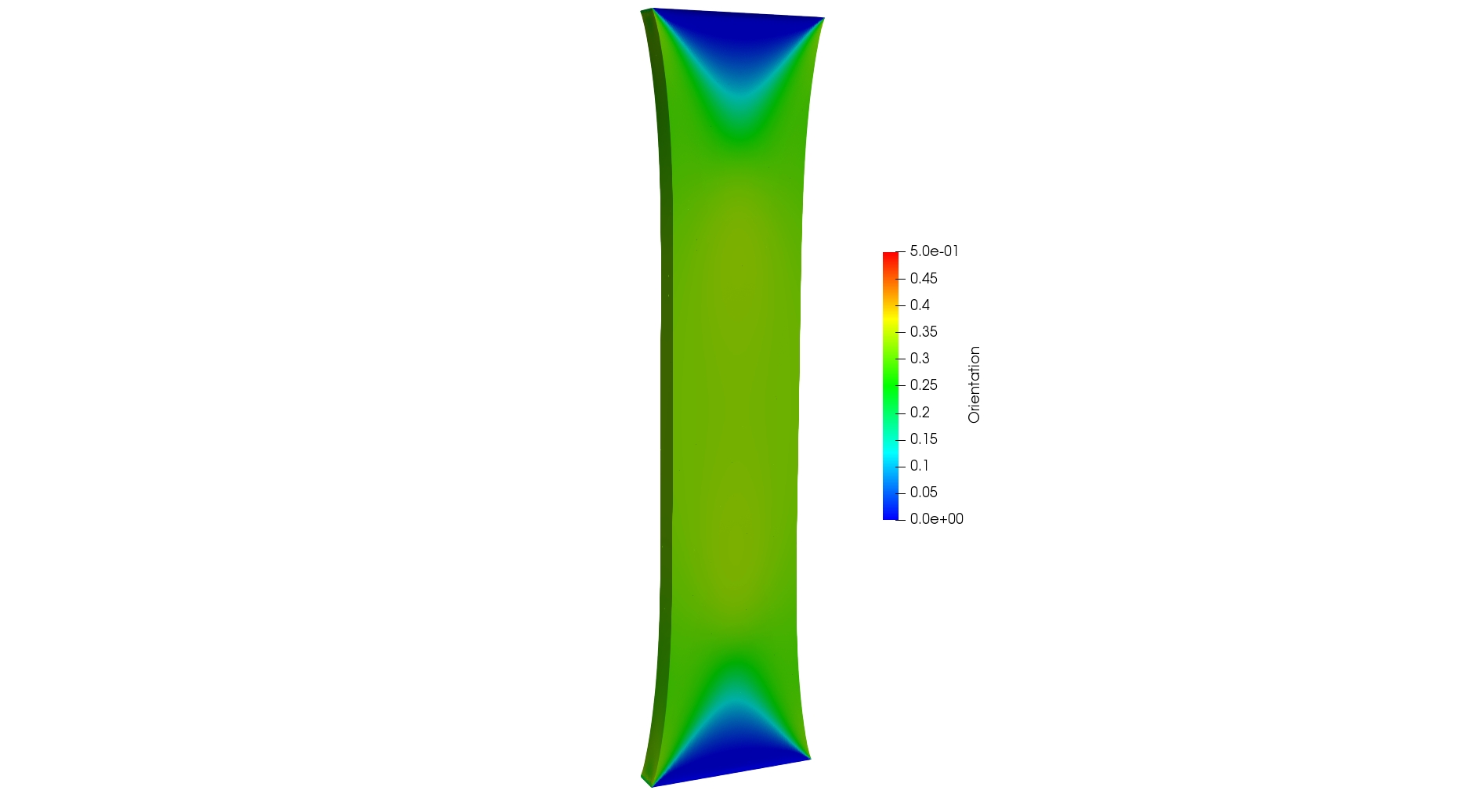} &
\includegraphics[angle=0, trim=800 30 700 30, clip=true, scale = 0.23]{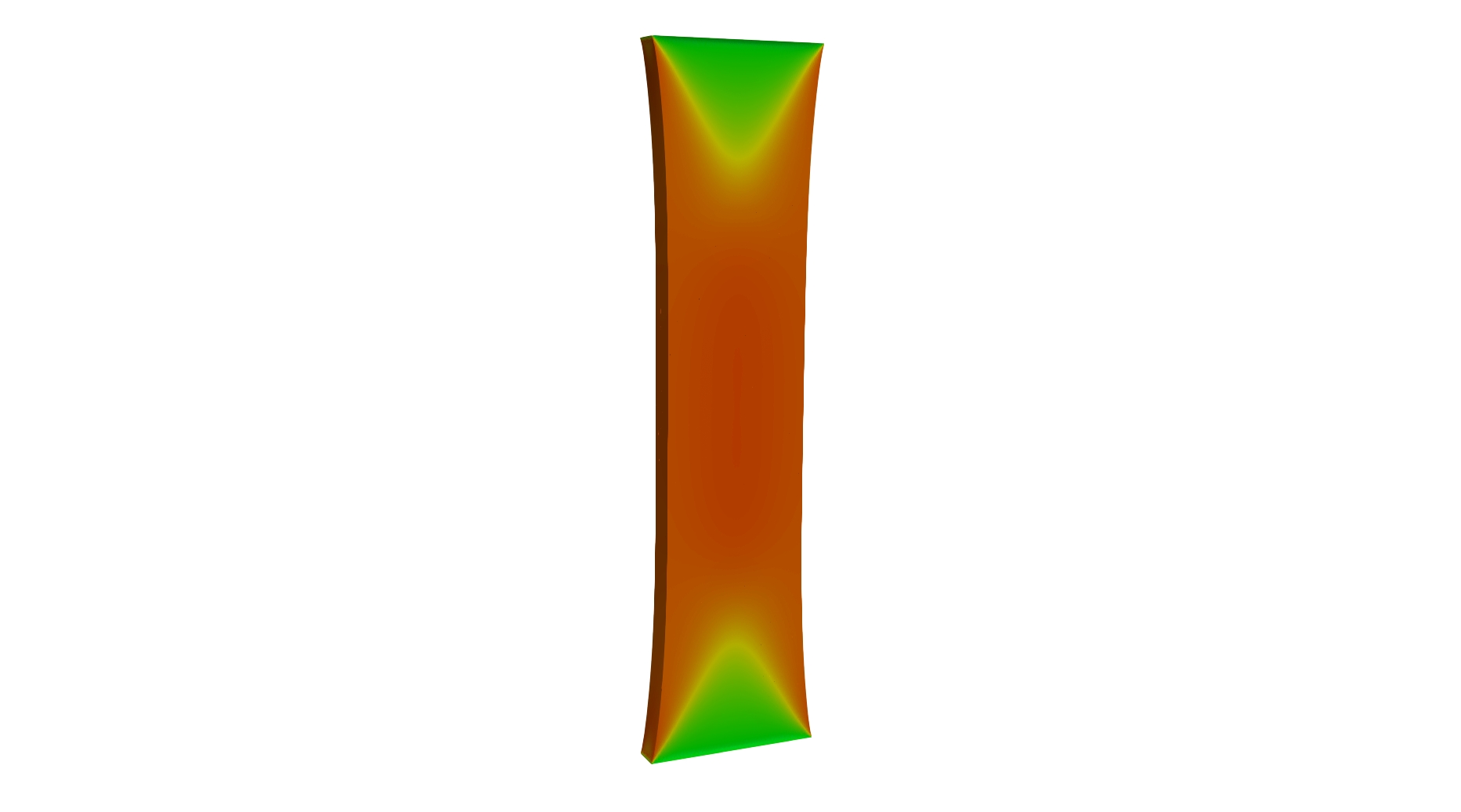} \\
(a) & (b) & (c) & (d)
\end{tabular}
\caption{Three-dimensional tensile test: Cauchy stress in the loading direction are plotted for the circumferential (a) and axial (b) specimens. The mean orientation of the collagen fibres in the current configuration are plotted for the circumferential (c) and axial (d) specimens.} 
\label{fig:tensile_test_stress}
\end{center}
\end{figure}

Before studying the solver performance, we perform a simulation with 3.5 million unstructured linear tetrahedral elements to examine the VMS formulation for this material model. In this study, the tensile test is performed in a dynamic approach. The loading force is applied as a linear function of time and reaches $2$ N in $100$ seconds. We set the density of the tissue as $1.0$ $\textup{g/cm}^3$. The tensile load-displacement curves for the circumferential and axial specimens with $k_d = 0.0$ and $0.226$ are plotted in Figure \ref{fig:tensile_test_setting} (b). We observe that before the fibres align along the loading direction, the groudmatrix provides the load carry capacity and the material response is very soft. When the fibres rotate to align with the loading direction, they take over the load burden, the material becomes stiffer, and the stiffness grows exponentially. For the axial specimen, the mean orientation of the fibres are closer to the loading direction, and hence it stiffens earlier than the circumferential specimen. Compared with the dispersed case, the specimen with perfectly aligned fibres (i.e., $k_d = 0.0$) needs a significant amount of rotation before they can carry load. In Figure \ref{fig:tensile_test_stress} (a) and (b), the Cauchy stresses in the tensile direction for the circumferential and axial specimens with $k_d = 0.226$ at the tensile load $1.0$ N are illustrated. The value of $\bm a_1 \cdot \bm C \bm a_2 / \|\bm F \bm a_1\| \|\bm F \bm a_2\|$ characterizes the current fibre alignment, and it is illustrated in Figure \ref{fig:tensile_test_stress} (c) and (d) for the circumferential and axial specimens. The maximum values in these specimens are $0.347$ and $0.459$, respectively. Correspondingly, the angles between the current mean fibre direction and the circumferential direction are $34.85^{\circ}$ and $31.34^\circ$, respectively. In the following discussion, the problem has been non-dimensionalized by the centimetre-gram-second units. Except the study performed in Section \ref{subsec:conv_vary_fiber_orientation_dispersion}, we adopt the axial specimen with the dispersion parameter $k_d = 0.226$ as the model problem for the study of the solver performance.

\subsubsection{Performance with varying inner solver accuracy}
\label{subsec:conv_rate_anisotropic_inner_solver}
In this test, we study the impact of the inner solver accuracy on the iterative solution algorithm. We fix the mesh size to be $1/400$ and the time step size to be $10^{-5}$. The simulation is performed with 8 CPUs. In this test, the settings of the linear solver are identical to the study performed in Section \ref{subsec:conv_rate_isotropic_compression}. The statistics of the solver are collected for the first time step of the simulation with varying values of $\delta^r_I$ (Table \ref{table:impact_inner_on_FGMRES_anisotropic}). We observe that using the inner solver may significantly improve the convergence rate of the linear solver. In both cases, the optimal performance in terms of time to solution is achieved by setting $ \delta^r_I = 10^2 \delta^r_S$.

\begin{table}[htbp]
\begin{center}
\tabcolsep=0.19cm
\renewcommand{\arraystretch}{1.2}
\begin{tabular}{P{3.0cm} P{1.0cm} P{3.0cm} P{0.8cm} P{0.8cm} P{1.0cm} P{1.0cm} P{1.0cm} }
\hline
& $\delta^r_I$ & CPU time (sec.) & $\hat{l}$ &  $n$ & $\bar{n}_A$ & $\bar{n}_S$ & $\bar{n}_I$   \\
\hline
$\delta^r_A = \delta^r_S = 10^{-10}$ & $10^{0}$ & $7.56 \times 10^1$ & 1 & 47 & 16.81 & 19.81 & -  \\
& $10^{-2}$ & $7.19 \times 10^1$ & 1 & 9 & 15.78 & 58.56 & 1.95  \\
& $10^{-4}$ & $6.52 \times 10^1$ & 1 & 5 & 15.30 & 58.80 & 3.75 \\
& $10^{-6}$ & $6.20 \times 10^1$ & 1 & 3 & 14.33 & 58.33 & 7.81  \\
& $10^{-8}$ & $5.83 \times 10^1$ & 1 & 2 & 13.75 & 59.50  & 11.85  \\
& $10^{-10}$ & $7.55 \times 10^1$ & 1 & 2 & 13.75 & 60.00 & 15.19  \\[0.5em]
$\delta^r_A = \delta^r_S = 10^{-6}$ & $10^{0}$ &  $4.38 \times 10^1$  & 1 & 47 & 7.83 & 12.02  & -  \\
& $10^{-2}$ & $4.20 \times 10^1$  & 1 & 9 & 8.17  & 34.22  & 1.97  \\
& $10^{-4}$ & $3.67 \times 10^1$ & 1 & 5 & 7.90 & 34.40 & 3.40 \\
& $10^{-6}$ & $4.73 \times 10^1$  & 1 & 4 & 7.50 & 35.75 & 7.37 \\
& $10^{-8}$ & $6.87 \times 10^1$ & 1 & 4 & 7.00 & 35.75 & 11.50  \\
& $10^{-10}$ & $8.23 \times 10^1$ & 1 & 4 & 7.00 & 35.75 & 15.47  \\
\hline
\end{tabular}
\end{center}
\caption{The impact of the accuracy of the inner solver on the performance of the linear solver. The CPU time is collected for the linear solver only; $\hat{l}$ represents the total number of nonlinear iterations; $n$ represents the total number of FGMRES iterations; $\bar{n}_A$ represents the averaged number of iterations for solving with $\boldsymbol{\mathrm A}$ in \eqref{eq:seg_sol_int_disp} and \eqref{eq:seg_sol_disp}; $\bar{n}_S$ represents the averaged number of iterations for solving \eqref{eq:seg_sol_pres}; $\bar{n}_I$ represents the averaged number of iterations for solving \eqref{eq:S_inner_A_eqn}.}
\label{table:impact_inner_on_FGMRES_anisotropic}
\end{table}

\subsubsection{Performance with varying intermediate solver accuracy}
\label{subsec:conv_rate_anisotropic_middle_solver}
\begin{figure}
	\begin{center}
	\begin{tabular}{cc}
\includegraphics[angle=0, trim=80 100 360 280, clip=true, scale = 0.04]{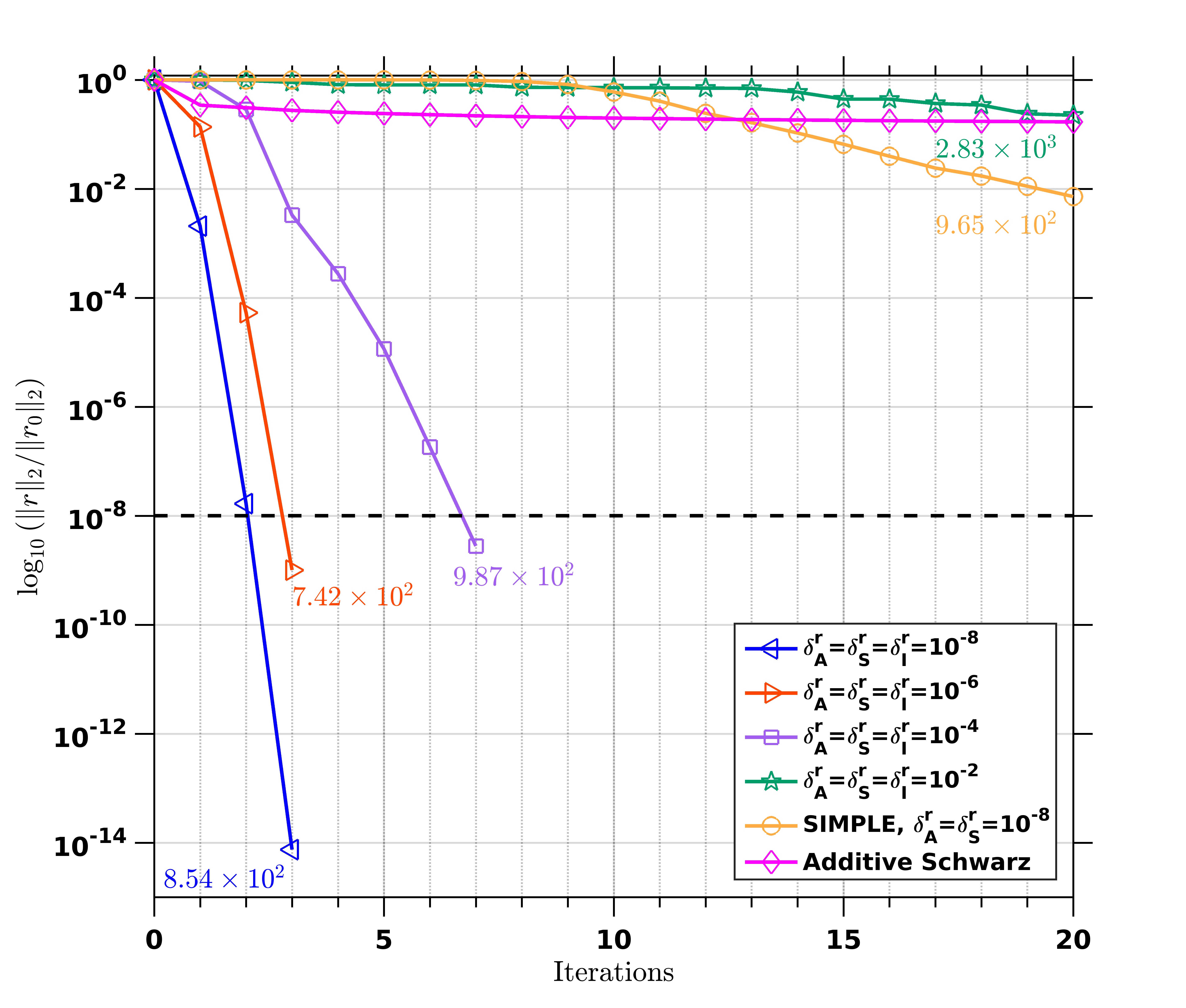} &
\includegraphics[angle=0, trim=80 100 360 280, clip=true, scale = 0.04]{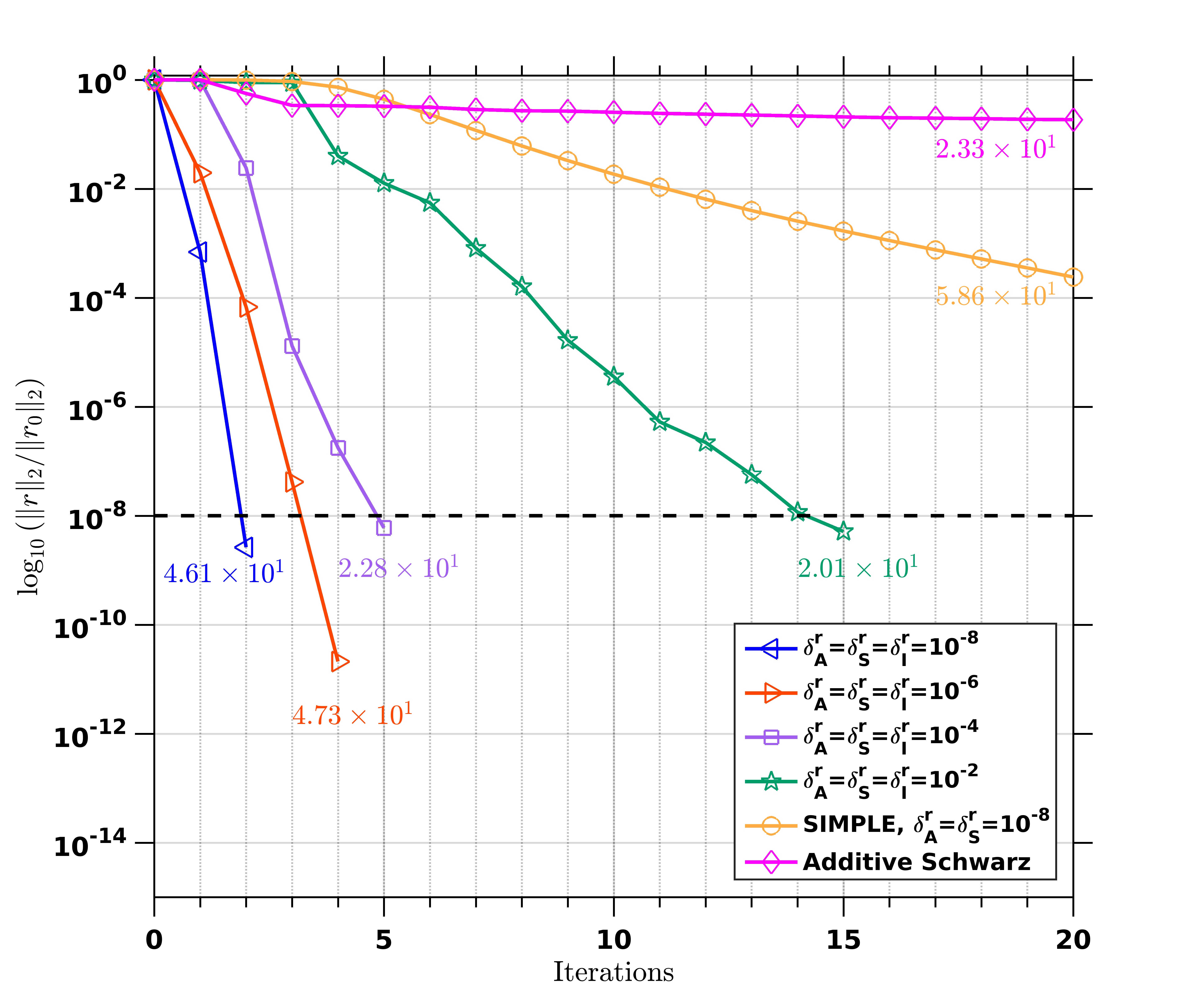}
\end{tabular}
\caption{Convergence history for $\Delta t = 10^{-1}$ (left) and $10^{-5}$ (right). The horizontal dashed black line indicates the prescribed stopping criterion for the relative error, which is $10^{-8}$ here. In the case of $\Delta t = 10^{-1}$, the block preconditioner with tolerance $10^{-2}$ converge in 90 iterations, the SIMPLE method converges in 45 steps, and the additive Schwarz method failed to converge. In the case of $\Delta t = 10^{-5}$, the SIMPLE method converges in 46 iterations, and the additive Schwarz method converges in 1070 iterations. The numbers indicate the averaged time per nonlinear iteration in seconds.} 
\label{fig:conv_history_ani_m10}
\end{center}
\end{figure}

We examine the solver performance for anisotropic hyperelastic materials with varying tolerances for the intermediate solvers. The mesh size is fixed to be $\Delta x = 1/400$, and the time step sizes are fixed to be $\Delta t = 10^{-1}$ and $10^{-5}$. The simulations are performed with 8 CPUs. We choose $\delta^r = 10^{-8}$ and vary the values of $\delta^r_A = \delta^r_S = \delta^r_I$ from $10^{-8}$ to $10^{-2}$. The SIMPLE preconditioner and the additive Schwarz preconditioner are also simulated for comparison. In the SIMPLE preconditioner, the block matrices $\boldsymbol{\mathrm A}$ and $\hat{\boldsymbol{\mathrm S}}$ are solved with $\delta^r_A = \delta^r_S = 10^{-8}$. The convergence history of the linear solver in the first nonlinear iteration is plotted in Figure \ref{fig:conv_history_ani_m10}. We observe that the nested block preconditioner performs robustly with a strict choice of the intermediate and inner solver tolerances. When the tolerances for the intermediate and inner solvers are loose ($10^{-2}$) and the time step is large ($\Delta t = 10^{-1}$), the convergence rate of the nested block preconditioner slows dramatically and is slower than the SIMPLE preconditioner. It should be emphasized that the SIMPLE preconditioner uses a very strict tolerance ($\delta^r_A = \delta^r_S = 10^{-8}$) here. We also note that the additive Schwarz preconditioner fails to converge to the prescribed tolerance in 10000 iterations when $\Delta t = 10^{-1}$.

\subsubsection{Performance with varying fibre orientations and dispersions}
\label{subsec:conv_vary_fiber_orientation_dispersion}
In this test, we examine the robustness of the solver with different collagen fibre orientations and dispersions. The structure of the arterial wall is described by the collagen fibre mean orientation $\varphi$ and the dispersion parameter $k_d$. We vary the value of $\varphi$ from $20^{\circ}$ to $80^{\circ}$, and the value of $k_d$ from $0.1$ to $0.3$. The rest material properties are kept the same as the ones used in the previous studies. The simulations are performed with $\Delta x = 1/100$ on $8$ CPUs. The time step size is $\Delta t = 10^{-1}$, and we simulate the problem up to $T=1.0$ to collect statistics of the solver performance. The stopping condition for the FGMRES iteration is $\delta^r = 10^{-8}$, and we choose $\delta^r_A = \delta^r_S = \delta^r_I = 10^{-6}$. The averaged number of iterations and the averaged CPU time for one nonlinear iteration is reported in Table \ref{table:material_dir_kd}.

\begin{table}[htbp]
\begin{center}
\tabcolsep=0.19cm
\renewcommand{\arraystretch}{1.2}
\begin{tabular}{P{2.5cm} P{4.1cm} P{4.1cm} P{4.1cm}}
\hline
$\bar{n}$ [$\bar{n}_A$, $\bar{n}_S$] ($\bar{T}_L$) & $k_d=0.1$ & $k_d=0.2$ & $k_d=0.3$ \\
 $\varphi = 20^\circ$ & 3.0 [214.1, 17.6] ($4.8\times 10^1$) & 3.0 [161.2, 18.0] ($3.4 \times 10^1$) & 3.0 [103.4, 19.7] ($2.1\times 10^1$)  \\
 $\varphi = 40^\circ$ & 3.0 [241.3, 17.7] ($5.8\times 10^1$) & 3.0 [176.8, 18.4] ($4.1 \times 10^1$) &  2.9 [105.7, 20.5] ($2.1\times 10^1$) \\
 $\varphi = 60^\circ$ & 2.8 [221.4, 17.8] ($4.6\times 10^1$) & 2.9 [169.1, 19.2] ($3.7 \times 10^1$)  & 2.9 [104.5, 20.9] ($2.1\times 10^1$)\\
 $\varphi = 80^\circ$ & 2.9 [220.8, 18.1] ($5.3\times 10^1$) & 3.0 [168.2, 19.8] ($4.1 \times 10^1$)  & 3.0 [103.0, 20.9] ($2.2\times 10^1$)\\
\hline
\end{tabular}
\end{center}
\caption{The performance of the nested block preconditioner with varying fibre orientations and dispersions.}
\label{table:material_dir_kd}
\end{table}

We observe that the outer solver converges in around three iterations regardless of the structural properties. In the intermediate level, the linear solver for $\boldsymbol{\mathrm S}$ is not sensitive to the two structural parameters; the linear solver for $\boldsymbol{\mathrm A}$ is affected by both parameters. The dispersion parameter $k_d$ has a significant impact on the performance of the solver associated with $\boldsymbol{\mathrm A}$. For the case of $k_d = 0.1$, the solver for $\boldsymbol{\mathrm A}$ requires slightly more than 200 iterations for convergence; for the  case of $k_d = 0.3$, the number of iterations drops to around 100.  As the dispersion parameter grows, there are more fibres providing stiffness. Thus, the trend of $\bar{n}_A$ is in agreement with the observations made in Section \ref{subsec:solver_different_material_para}.

\subsubsection{Parallel performance}
\label{subsec:anisotropic_parallel_performance}
We compare the performance of different preconditioners by performing a weak scaling test. The tolerance for the linear solver is set to be $\delta^r = 10^{-3}$. In the nested block preconditioner, we set $\delta^r_A = \delta^r_S = \delta^r_I = 10^{-3}$, and we use $\delta^r_A = \delta^r_S = 10^{-3}$ for the SIMPLE preconditioner. The computational mesh is progressively refined and each CPU is assigned with approximately $6.0 \times 10^4$ equations. We simulate the problem with two different time step sizes: $\Delta t = 10^{-1}$ and $10^{-5}$. The statistics of the solver performance are collected for five time steps, and the results are reported in Table \ref{table:ani_weak_scaling}. The number of iterations at the intermediate level shows a similar trend to the isotropic case studied in Section \ref{subsec:parallel_performance}. The difference is that, for the anisotropic material, the solver for $\boldsymbol{\mathrm A}$ requires more iterations to converge when the time step size is large. The degradation of the AMG preconditioner for anisotropic problems is known, and using a higher complexity coarsening, like the Falgout method, will improve the performance \cite{Henson2002}. Notably, for large time steps, the additive Schwarz preconditioner just cannot deliver converged solutions within 10000 iterations, regardless of the spatial mesh size. Examining the results, the proposed nested block preconditioner gives the most robust and efficient performance for most of the cases considered.

\begin{table}
\begin{center}
\tabcolsep=0.19cm
\renewcommand{\arraystretch}{1.2}
\begin{tabular}{@{\extracolsep{4pt}}P{1.0cm} P{1.0cm} P{0.6cm} P{1.0cm} P{1.0cm} P{1.0cm} P{1.0cm} P{1.0cm} P{1.5cm} P{1.0cm}@{}}
\hline
\multirow{2}{*}{$\frac{1}{\Delta x}$} & \multirow{2}{*}{Proc.} & \multicolumn{4}{c}{$\hat{\mathcal P}_{SCR}$} & \multicolumn{2}{c}{SIMPLE} & \multicolumn{2}{c}{Additive Schwarz}\\
\cline{3-6} \cline{7-8}\cline{9-10}
& & $\bar{n}$ & $\bar{n}_A$ & $\bar{n}_S$ & $\bar{T}_L$ & $\bar{n}$ & $\bar{T}_L$ & $\bar{n}$ & $\bar{T}_L$
\\
\hline
\multicolumn{10}{l}{$\Delta t = 10^{-1}$} \\
$200$ & 8 & 6.2 & 207.1 & 10.3 & 579.3  & 54.1 & 888.5 & NC & NC \\
$400$ & 64 & 7.8 & 331.3  & 10.5 & 2062.7 & 102.5 &  6262.8 & NC & NC \\
$600$ & 216 & 10.3  & 389.4 & 11.1 & 3204.1 & 140.5 & 8051.3 & NC & NC \\ 
\hline
\multicolumn{10}{l}{$\Delta t = 10^{-5}$} \\
$200$ & 8 & 4.0 & 2.3 & 15.5 & 9.7 & 22.6 & 9.4 & 485.6 & 10.8 \\  
$400$ & 64 & 4.1 & 3.1 & 18.6 & 29.2 & 45.3 & 53.3 & 986.3 & 47.76 \\ 
$600$ & 216 & 5.9 & 3.9  & 20.3 & 119.8 & 71.3 & 202.5 & 1453.0  & 431.6 \\
\hline
\end{tabular}
\end{center}
\caption{Comparison of the averaged iteration counts and CPU time in seconds for the nested block preconditioner $\hat{\mathcal P}_{SCR}$, the SIMPLE preconditioner, and the additive Schwarz preconditioner. NC stands for no convergence. For the $\Delta t = 10^{-1}$ case, the additive Schwarz preconditioner failed to achieve convergence in 10000 iterations.}
\label{table:ani_weak_scaling}
\end{table}

\section{Conclusions}
\label{sec:conclusion}
In this work, we designed a preconditioning technique based the novel hyper-elastodynamics formulation \cite{Liu2018}. This preconditioning technique is based on a series of block factorizations in the Newton-Raphson solution procedure \cite{Liu2018,Scovazzi2016,Rossi2016} and is inspired from the preconditioning techniques developed in the CFD community \cite{Bank1990,Baggag2004,Manguoglu2008,Manguoglu2009}. It uses the Schur complement reduction with relaxed tolerances as the preconditioner inside a Krylov subspace method. This strategy enjoys the merits of both the SCR approach and the fully coupled approaches. It shows better robustness and efficiency in comparison with the SIMPLE and the additive Schwarz preconditioners. Tuning the intermediate and the inner solvers allows the user to adjust the nested algorithm for specific problems to attain a balance between robustness and efficiency. In this work, to make the presentation coherent, we adopted the same solver at the intermediate and the inner levels. In practice, one is advised to flexibly apply the most efficient solver at the inner level. For example, one may symmetrize the matrix in \eqref{eq:S_inner_A_eqn} \cite{Manguoglu2009} and use the conjugate gradient method as the inner solver. In our experience, this will further reduce the computational cost. In all, the methodology developed in this work provides a sound basis for the design of effective preconditioning techniques for hyper-elastodynamics.

There are several promising directions for future work. (1) Improvements will be made to design a better preconditioner for the Schur complement. It is tempting to consider using the sparse approximate inverse method to construct this preconditioner \cite{Chow2000}. (2)  This preconditioning technique will be extended to inelastic calculations \cite{Zeng2017,Abboud2018} as well as FSI problems \cite{Liu2018}.

\section*{Acknowledgements}
This work is supported by the National Institutes of Health under the award numbers 1R01HL121754 and 1R01HL123689, the National Science Foundation (NSF) CAREER award OCI-1150184, and computational resources from the Extreme Science and Engineering Discovery Environment (XSEDE) supported by the NSF grant ACI-1053575. The authors acknowledge TACC at the University of Texas at Austin for providing computing resources that have contributed to the research results reported within this paper.

\appendix
\section{Consistent linearization}
\label{app:consistent_linearization}
We report the explicit formulas of the residual vectors and tangent matrices used in the Newton-Raphson solution procedure at the iteration step $l$. For notational simplicity, the subscript $(l)$ is neglected in the following discussion.
\begin{align}
\boldsymbol{\mathrm R}_{p} =& \left[ \mathrm R_{p,A} \right], \displaybreak[2] \\
\mathrm R_{p,A} =& \int_{\Omega_{\bm X}}  J N_A \left( \beta \dot{p} + v_{i,i} \right) d\Omega_{\bm X} + \sum _{e}\int_{\Omega_{\bm X}^e} \tau^e_{M} N_{A,i} \left( \rho J \dot{v}_i - \tilde{\bm P}_{iJ,J} + J p_{,i} - \rho J b_i \right)  d\Omega_{\bm X}, \displaybreak[2] \\
\boldsymbol{\mathrm R}_{m} =& \left[ \mathrm R_{m,A}^i \right], \displaybreak[2] \\
\mathrm R_{m,A}^i =& \int_{\Omega_{\bm X}} N_A J \rho \dot{v}_i + N_{A,I} \tilde{\bm P}_{iI} - N_{A,i} J p - N_A J\rho b_i d\Omega_{\bm X} - \int_{\Gamma^H_{\bm X}} N_A H_i d\Gamma_{\bm X}.
\end{align}
In the above, we used the following notation conventions,
\begin{align*}
N_{A,I} := \frac{\partial N_A}{\partial X_I}, \quad N_{A,i} := \frac{\partial N_A}{\partial x_i} = \frac{\partial N_A}{\partial X_I} \frac{\partial X_I}{\partial x_i} = \frac{\partial N_A}{\partial X_I} F^{-1}_{Ii}, \quad p_{,i} = p_{,I} F^{-1}_{Ii}, \quad \tilde{\bm P}_{iI} := J \bm \sigma^{dev}_{ij} \bm F^{-1}_{Ij}.
\end{align*}
Note that $\rho = \rho(p)$ and $\beta = \beta(p)$ are given by the constitutive relations, and $H_i := h_i \circ \bm \varphi_t$.
\begin{align}
\boldsymbol{\mathrm A} =& \left[ \mathrm A_{AB}^{ij} \right], \displaybreak[2] \\
\mathrm A_{AB}^{ij} =& \alpha_m \int_{\Omega_{\bm X}} J \rho N_A N_B d\Omega_{\bm X} \delta_{ij} + \frac{\left(\alpha_f \gamma \Delta t_n \right)^2}{\alpha_m} \int_{\Omega_{\bm X}} N_{A,I} \left( \tilde{\bm S}_{IJ}\delta_{ij} + \mathbb A^{iso}_{iIjJ} \right) N_{B,J} d\Omega_{\bm X}  \nonumber \displaybreak[2] \\
& + \frac{\left(\alpha_f \gamma \Delta t_n \right)^2}{\alpha_m} \int_{\Omega_{\bm X}} J p N_{A,I} \left( F^{-1}_{Ij} F^{-1}_{Ji} - F^{-1}_{Ii} F^{-1}_{Jj} \right) N_{B,J} d\Omega_{\bm X}  \nonumber \displaybreak[2] \\
& + \frac{\left(\alpha_f \gamma \Delta t_n \right)^2}{\alpha_m} \int_{\Omega_{\bm X}} N_A \rho J \left( \dot{v}_i - b_i \right) N_{B,j} d\Omega_{\bm X}, \nonumber \displaybreak[2] \\
\boldsymbol{\mathrm B} =& \left[ \mathrm B_{AB}^{i} \right], \displaybreak[2] \\
\mathrm B_{AB}^{i} =& \alpha_f \gamma \Delta t \int_{\Omega_{\bm X}} \rho_{,p}J (\dot{v}_i -b_i) N_A N_B - J N_{A,i} N_B d\Omega_{\bm X}, \nonumber \displaybreak[2] \\
\boldsymbol{\mathrm C} =& \left[ \mathrm C_{AB}^{j} \right], \displaybreak[2] \\
\mathrm C_{AB}^{j} =& \alpha_m \sum_e \int_{\Omega_{\bm X}^e} \tau^e_M \rho J N_{A,j} N_B d\Omega_{\bm X} + \alpha_f \gamma \Delta t \int_{\Omega_{\bm X}}  J N_A N_{B,j} d\Omega_{\bm X} \nonumber \displaybreak[2] \\
& + \frac{\left(\alpha_f \gamma \Delta t_n \right)^2}{\alpha_m} \int_{\Omega_{\bm X}} J \beta \dot{p} N_A N_{B,j} + J N_A \left( v_{i,i}N_{B,j} - v_{i,j} N_{B,i} \right) d\Omega_{\bm X} \nonumber \displaybreak[2] \\
& + \frac{\left(\alpha_f \gamma \Delta t_n \right)^2}{\alpha_m} \sum_e \int_{\Omega_{\bm X}^e} \tau^e_M N_{A,j} N_{B,i} \left( \rho J \dot{v}_i - \tilde{\bm P}_{iI,I} + J p_{,i} - \rho J b_i \right)d\Omega_{\bm X} \nonumber \displaybreak[2] \\
& + \frac{\left(\alpha_f \gamma \Delta t_n \right)^2}{\alpha_m} \sum_e \int_{\Omega_{\bm X}^e} \tau^e_M J N_{A,i} \left( p_{,i}N_{B,j} - p_{,j} N_{B,i} \right) d\Omega_{\bm X} \nonumber \displaybreak[2] \\
& + \frac{\left(\alpha_f \gamma \Delta t_n \right)^2}{\alpha_m} \sum_e \int_{\Omega_{\bm X}^e} \tau^e_M N_{A,i} N_{B,j} \rho J \left( \dot{v}_i - b_i \right) d\Omega_{\bm X} \nonumber \displaybreak[2] \\
& - \frac{\left(\alpha_f \gamma \Delta t_n \right)^2}{\alpha_m} \sum_e \int_{\Omega_{\bm X}^e} \tau^e_M N_{A,i} \left( \tilde{\bm S}_{MN} \delta_{ij} + \mathbb A^{iso}_{iMjN}  \right) N_{B,MN} d\Omega_{\bm X}, \displaybreak[2] \\
\boldsymbol{\mathrm D} =& \left[ \mathrm D_{AB} \right], \displaybreak[2] \\
\mathrm D_{AB} =& \alpha_m \int_{\Omega_{\bm X}} J \beta N_A N_B d\Omega_{\bm X} + \alpha_f \gamma \Delta t \int_{\Omega_{\bm X}} J \beta_{,p}\dot{p} N_A N_B d\Omega_{\bm X} \nonumber \displaybreak[2] \\
& + \alpha_f \gamma \Delta t \sum_e \int_{\Omega_{\bm X}^e} J \tau^e_M \left( N_{A,i} N_{B,i} + \rho_{,p} \left(\dot{v}_i - b_i\right) N_{A,i}N_B \right)d\Omega_{\bm X}.
\end{align}
In $\mathrm A^{ij}_{AB}$ and $\mathrm C^{j}_{AB}$, we used the following notation,
\begin{align*}
\mathbb A^{iso}_{iIjJ} := \frac{\partial G_{iso}}{\partial F_{iI} \partial F_{jJ}}.
\end{align*}

\bibliographystyle{plain}  

\bibliography{Solver_Solid}

\end{document}